\documentclass[a4paper,12pt,twoside]{article}
\usepackage[latin1]{inputenc}                 
\usepackage{amsfonts,amssymb,amsmath,exscale} 
\usepackage[dvips]{graphicx,color,psfrag}     
\usepackage{fancyhdr}                         
\usepackage{caption}                          
\usepackage{rotating}
\usepackage{defcml}                           
\usepackage[numbers]{natbib}                  
\Floatboxname{Box}
\usepackage{verbatim}                         

\usepackage{amsmath,amsthm,amsfonts,amssymb,amscd}
\usepackage{subfigure}
\usepackage{graphicx,bm,color}
\usepackage{algorithm}
\usepackage{algpseudocode}
\usepackage{stmaryrd}
\usepackage{booktabs,ctable,multirow,longtable} 
\usepackage[latin1]{inputenc}                 
\usepackage{amsfonts,amssymb,amsmath,exscale} 
\usepackage[dvips]{graphicx,color,psfrag}     
\usepackage{fancyhdr}                         
\usepackage{caption}                          
\usepackage{rotating}
\usepackage{defcml}                            
\usepackage[numbers]{natbib}                  
\Floatboxname{Box}

\definecolor{gray}{gray}{0.6}

\newtheorem{Remark}{Remark}[section]

\newtheorem{form}{Formulation}[section]

\def\R{\mathbb R}

\fancypagestyle{plain}{%
	\fancyhead{}%
	\fancyfoot[c]{\sffamily\thepage}%
}
\makeatletter
\def\cleardoublepage{\clearpage\if@twoside \ifodd\c@page\else
	\hbox{}
	\vspace*{\fill}
	\thispagestyle{empty}
	\newpage
	\if@twocolumn\hbox{}\newpage\fi\fi\fi}

\usepackage{scalerel,stackengine}
\stackMath
\newcommand\reallywidecheck[1]{%
	\savestack{\tmpbox}{\stretchto{%
			\scaleto{%
				\scalerel*[\widthof{\ensuremath{#1}}]{\kern-.6pt\bigwedge\kern-.6pt}%
				{\rule[-\textheight/2]{1ex}{\textheight}}
			}{\textheight}%
		}{0.5ex}}%
	\stackon[1pt]{#1}{\scalebox{-1}{\tmpbox}}%
}

\begin{document}
\Titel{
Global-local techniques for hydraulic fracture
      }
\Autor{F. Aldakheel, N. Noii, T. Wick, M. Wheeler, P. Wriggers}
\Report{02--I--17}
\Journal{

}
%



\thispagestyle{empty}

\ce{\bf\large A global-local approach for hydraulic phase-field fracture}
\vskip .1 in
\ce{\bf\large in poroelastic media}
\vskip .35in

\ce{
Fadi Aldakheel\(^{a}\), Nima Noii\(^{b,}\)\footnote{Corresponding author.\\[1mm] 
	E-mail addresses:
	aldakheel@ikm.uni-hannover.de (F. Aldakheel); noii@ifam.uni-hannover.de (N. Noii); thomas.wick@ifam.uni-hannover.de (T. Wick); wriggers@ikm.uni-hannover.de (P. Wriggers).
}, Thomas Wick\(^{b,c}\), Peter Wriggers\(^{a,c}\)} \vskip .25in

\ce{\(^a\) Institute of Continuum Mechanics} \ce{Leibniz Universit\"at Hannover, Appelstrasse 11, 30167 Hannover, Germany}\vskip .25in

\ce{\(^b\) Institute of Applied Mathematics} \ce{Leibniz Universit\"at Hannover, Welfengarten 1, 30167 Hannover, Germany} \vskip .25in

\ce{\(^c\) Cluster of Excellence PhoenixD (Photonics, Optics, and
	Engineering - Innovation} \ce{Across Disciplines) Leibniz Universit\"at Hannover, Germany}\vskip .25in

\begin{Abstract}
In this work, phase-field modeling of hydraulic fractures in porous media is extended towards a global-local approach. Therein, the failure behavior is solely analyzed in a (small) local domain. In the surrounding medium, a simplified and linearized system of equations is solved. Both domains are coupled by Robin-type interface conditions. The fracture(s) inside the local domain are allowed to propagate and consequently both subdomains change within time. Here, a predictor-corrector strategy is adopted in which the local domain is dynamically adjusted to 
the current fracture pattern. The resulting framework is algorithmically described in detail and substantiated with some numerical tests.
	
\end{Abstract} 

\sectpa[Section1]{Introduction}

In recent years, several pressurized 
\cite{BourChuYo12,MiWheWi19,MiWheWi15b,WickLagrange2014,Wick15Adapt,HeiWi18_pamm,singh2018finite,NoiiWick2019} 
and fluid-filled  
\cite{MiWheWi14,WiSiWhe15,LeeWheWi16,Miehe2015186,MieheMauthe2015,ehlers17,heider2019phase,HEIDER2018116,LeeMinWhe2018,wang2017unified,LeeMiWheWi16,Cajuhi2017,LeeWheWiSri17,CHUKWUDOZIE2019957,Wilson2016264,aldakheelspp2020} 
phase-field fracture formulations 
have been proposed in the literature.
These studies range from 
modeling of pressurized and fluid-filled fractures, 
mathematical analysis, numerical modeling and simulations 
up to high-performance parallel computations. Recently various extensions 
towards multiphysics phase-field fracture in porous media have been proposed in which 
various phenomena couple as for instance proppant \cite{LeeMiWheWi16},
two-phase flow formulations \cite{LeeMiWheWi18} or given temperature variations \cite{NoiiWick2019}.
All these examples demonstrate the potential of phase-field for crack propagation.

Phase-field fracture is a regularized approach, which has advantages and shortcomings. The first advantage is a continuum description based on
first physical principles to determine the unknown crack path
\cite{FraMar98,BourFraMar00,miehe+hofacker+schaenzel+aldakheel15} and the computation 
of curvilinear and complex crack patterns. The model allows for
nucleation, branching, merging and post-processing of certain quantities
such that stress intensity factors become redundant. Therefore, easy handling
of fracture networks in possibly and highly heterogeneous media can be treated.
The formulation being described in a variational framework allows 
finite element discretizations and corresponding analyses. The mathematical model
permits any dimension, thus phase-field fracture applies conveniently to three-dimensional
simulations. On the energy level, the formulation is non-convex 
constituting a challenge for both theory and design of numerical algorithms.  
A second challenge is the computational cost. Various solutions 
have been proposed so far; namely staggered approaches (alternating minimization)
\cite{Bour07,BuOrSue10,BuOrSue13}, 
stabilized staggered techniques \cite{BrunEtAl19}, 
quasi-monolithic approaches \cite{Wick15Adapt} (possibly 
with sub-iterations \cite{MaWi19}), or fully monolithic approaches 
\cite{Gerasimov16,Wi17_SISC,Wi17_CMAME}.
Adaptive mesh refinement was proposed to reduce the computational costs 
\cite{BuOrSue10,Wick15Adapt,ArFoMiPe15,Wi16_dwr_pff}. A related technique that has the potential 
to treat large-scale problems is a global-local technique proposed 
in \cite{NoiiGL18}. Recently this was extended to a 
framework in which the local domain is dynamically updated according 
to the propagating fracture path \cite{NoiiAldakheelWickWriggers2019}.
The need for such framework can be found in {multiscale porous media 
	applications \cite{DaWhe18,DANA20181}} or in which a localized fracture occurs in a (big) reservoir
\cite{WiSiWhe15,GiWheAlDa19}. 

The last two references are the motivation for the present work. Here,
we extend the adaptive global-local phase-field fracture approach 
\cite{NoiiAldakheelWickWriggers2019} 
to porous media applications with hydraulic fractures.
We first extend our model towards large strain formulations, in line with \cite{aldakheel+mauthe+miehe14,aldakheel16,dittmann+fadi+etal18,miehe2017phase}. Previous studies 
only concentrated on small strain applications. Then, the coupled 
multiphysics fracture framework is carefully derived. Both subdomains 
will be coupled via Robin-type interface conditions, see \cite{NoiiAldakheelWickWriggers2019}. 
This leads to Lagrange multiplier formulations that are demanding from a mathematical 
point of view as well as in the implementation, see for example \cite{Wohlmuth2000,Wohlmuth2003,Wo11,seitz2016isogeometric}. A future rigorous numerical analysis of our global-local approach can be achieved with similar methodologies
as used in  \cite{GePeWheWi09,GiPeWheWi11}.
In particular, our formulation can deal with non-matching grids at the interface,
which is very interesting for cases {towards practical field problems} 
as mentioned in \cite{WiSiWhe15,GiWheAlDa19}
in which possibly various programming codes must be coupled. 
On the fine-scale level all (nonlinear) equations are solved. On the global level, 
only coarse representations of the pressure and crack phase-field are considered. 
As mentioned in the previous descriptions and references such multiphysics fracture formulations 
are challenging from a mathematical and numerical point of view. For these reasons, we concentrate in this paper 
on careful algorithmic descriptions including supporting numerical simulations. 
Here, our emphasis is on results, demonstrating the computational convergence 
properties of our proposed numerical schemes. A rigorous numerical analysis 
must be left for future work.

The outline of this paper is as follows: In Section \ref{Section2}, 
the governing equations are described. Then, in Section \ref{Section3},
the extension to a global-local formulation 
for pressurized fractures is derived. Therein, the Robin-type interface 
conditions are carefully discussed. This is followed by the final 
global-local algorithm. Afterward, we also discuss the dynamic choice 
of the local domains with the help of a predictor-corrector scheme.
In Section \ref{Section5} some numerical tests are carried out in order 
to substantiate our algorithmic developments.
\sectpa[Section2]{Phase-field modeling of hydraulic fracture}
This section outlines a theory of hydraulic phase-field fracture in poroelastic media undergoing finite strains. The constitutive formulations are based on three governing equations for the mechanical deformation, fluid pressure and the crack phase-field. Strong and weak formulations of the mutli-physics problem are introduced. Furthermore the framework is algorithmically described, resulting in the so-called {\it single-scale domain} formulations. 

\sectpb[Section2.1]{Governing equations}
\begin{figure}[!b]
	\centering
	{\includegraphics[clip,trim=5.5cm 7cm 2cm 9cm, width=16.5cm]{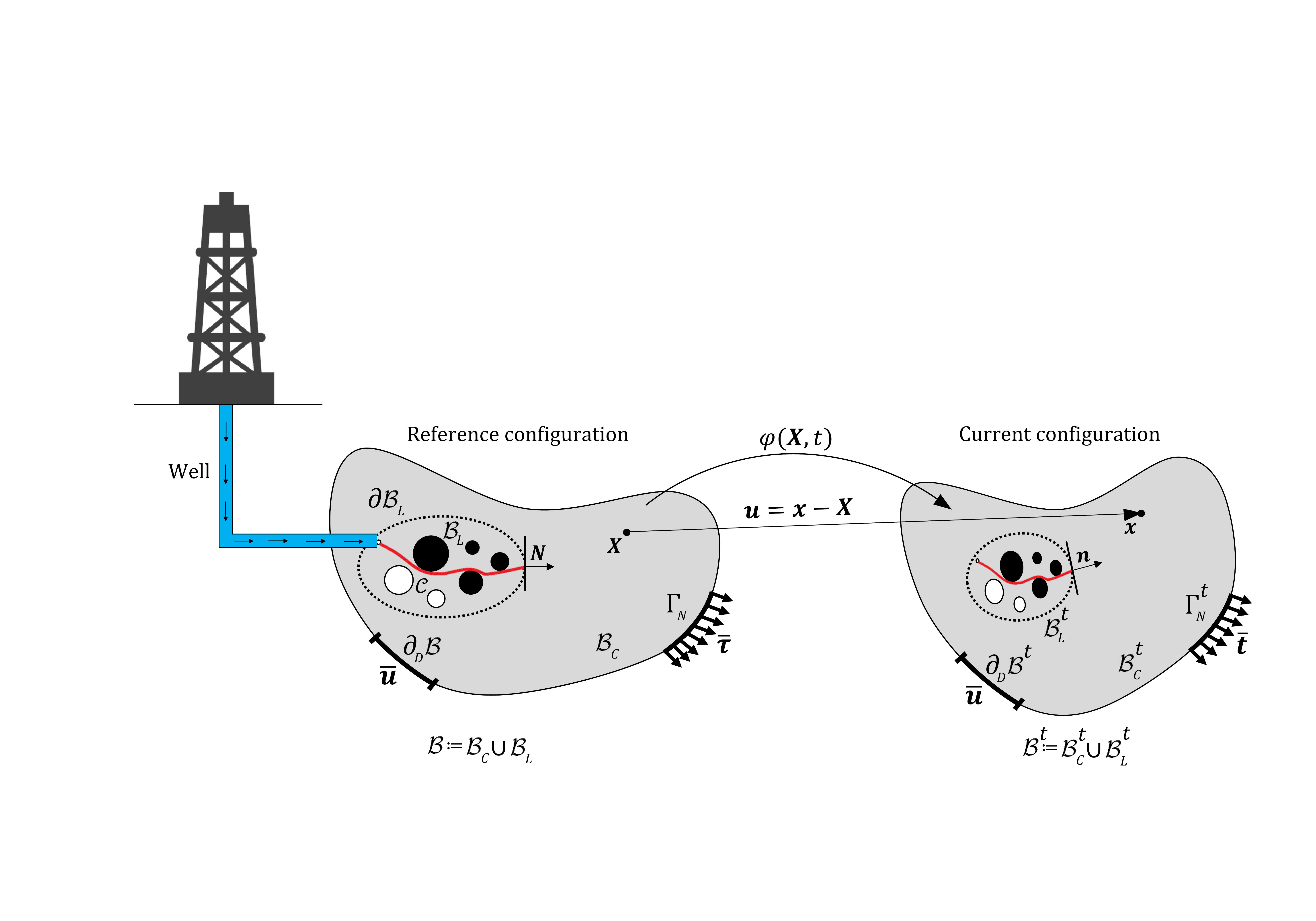}}  
	\caption{ Setup of the notation: the unbroken domain is denoted by
		$\calB_C$ and $\mathcal{C}$ is the crack phase-field. The smeared crack
		phase-field is approximated by the domain $\calB_L$. The whole domain is defined as a close subset as ${\calB:={\calB}_C\cup{\calB}_L}$. The fracture boundary is $\partial \calB_L$ and the outer boundary of the domain is $\partial
		\calB$. Blue color refer to the injected fluid through the well drilling and injection.}
	\label{Figure1}
\end{figure}
Consider $\calB\in{\calR}^{\delta}$ to be a material body (solid in the reference configuration) and denote $\partial\calB$ as its boundary with dimension $\delta = 2,3$ in space and time $t\in \calT = [0,T]$. We assume Dirichlet boundary conditions on $\partial_D\calB $ and Neumann conditions on $\partial_N \calB := \Gamma_N \cup \mathcal{C}$, where $\Gamma_N$  denotes the outer domain boundary. The lower dimensional curved surface $\calC\in \calR^{\delta-1}$ is the crack boundary, as illustrated in Fig. \ref{Figure1}. 

The boundary-value-problem BVP for the coupled problem of fluid-saturated porous media at fracture is a coupled three-field problem. It is characterized at material points $\BX\in\mathcal{B}$ by the deformation map $\Bvarphi(\BX,t)$ of the solid, the fluid pressure field $p(\BX,t)$ and the crack phase-field $d(\BX,t)$ defined as
\begin{equation}
\Bvarphi: 
\left\{
\begin{array}{ll}
\calB \times \calT \rightarrow \calR^{\delta} \\
(\BX, t)  \mapsto \Bx = \Bvarphi(\BX,t)
\end{array}
\right.
, \
p:
\left\{
\begin{array}{ll}
\calB \times \calT \rightarrow \calR \\
(\BX, t)  \mapsto p(\BX,t)
\end{array}
\right.
, \
d: 
\left\{
\begin{array}{ll}
\calB \times \calT \rightarrow [0,1] \\
(\BX, t)  \mapsto d(\BX,t)
\end{array}
\right.,
\label{phi-p-d-fields}
\end{equation}
with $\dot{d} \ge 0$. The position of a material point in the deformed configuration is depicted as $\Bx= \BX + \Bu(\BX,t)$ where $\Bu(\BX,t)$ is the displacement field. The crack phase-field $d(\BX,t)=0$ (light gray color) and $d(\BX,t)=1$ (red color) refer to the unbroken and fully fractured state of the material respectively, as visualized in Fig. \ref{Figure1}. The fracture surface $\mathcal{C}$ is approximated in $\calB_L\subset\calB$ the so-called \textit{local domain}. The intact region with no fracture is denoted
as \textit{complementary domain} $\calB_C:=\calB \backslash \calB_L\subset\calB$, such that ${\calB}_C\cup{\calB}_L=:\calB$ 
and ${\calB_C}\cap{\calB_L}=\varnothing$. We note that $\calB_L$, the domain in which the smeared crack phase-field is approximated, and its boundary $\partial \calB_L$ depend on the choice of the phase-field regularization parameter $l>0$.

The material deformation gradient of the solid is defined by $\BF :=
\nabla\Bvarphi(\BX,t)=\text{Grad}\Bvarphi$ with the Jacobian $J \!:=\!\mbox{det}[\BF]  > 0$ and the right Cauchy-Green tensor $\BC = \BF^T \BF$. The solid is loaded by prescribed deformations and external traction on the
boundary, defined by time-dependent Dirichlet- and Neumann conditions
\begin{equation}
\Bvarphi = \bar\Bvarphi(\BX,t) \ \textrm{on}\ \partial_D\calB
\AND
\BP\BN = \bar\Btau(\BX,t)\ \textrm{on}\ \partial_N\calB,
\label{mech-bcs}
\end{equation}
where $\BN$ is the outward unit normal vector and $\bar{\Btau}$ is the prescribed traction vector at the surface $\partial\calB$ of the undeformed configuration. The first Piola-Kirchoff stress tensor $\BP$ is the thermodynamic dual to $\BF$. 

The solid has to satisfy the equation of equilibrium, representing the {\it first} partial differential equation PDE for the coupled problem as
\begin{equation}
\boxed{\Div\,\BP + \overline\Bb = \Bzero}
\label{equil:defo}
\end{equation}
where dynamic effects are neglected and $\overline\Bb$ is the given body force.

For the constitutive modeling of poromechanics, we proceed with a {\it biphasic} fully saturated porous material, consisting of a pore fluid and a solid matrix material. A local volume element $dV$ in the undeformed reference configuration can be decomposed into a fluid $dV_{F}$ and a solid $dV_{S}$ part. Thereby the volume fraction can be defined via $n_{\alpha}:=dV_{\alpha}/dV$, where $\alpha=\{S,F\}$. The \textit{saturation condition} for the case of a fully saturated porous medium is given by
\begin{equation}
\sum_{\alpha}n_{\alpha} = n_{F}+n_{S} = 1\, ,
\end{equation}
where $n_{F}(\BX,t)$ represents the porosity, i.e. the volume occupied by the fluid is {\it same as} the pore volume. Note that in the {\it fracture region} where
\begin{equation}
d = 1 \quad \mbox{yields} \quad n_S =0 \AND n_F = 1 \ .
\end{equation}
The volume fraction in porous media relates the real density (material, effective, intrinsic) $\rho_{\alpha R}$ to the partial density $\rho_{\alpha}$ by
\begin{equation}
\rho_{\alpha}=n_{\alpha}\;\rho_{\alpha R} \quad \text{with} \quad \rho_{\alpha R}:=dm_{\alpha}/dV_{\alpha} \quad \text{and} \quad \rho_{\alpha}:=dm_{\alpha}/dV \; ,
\end{equation}
where $dm_{\alpha}$ is the mass of the phase $\alpha$. Thus, the  overall density can be expressed as $\rho=\sum_{\alpha}n_{\alpha} \;\rho_{\alpha R}$ .  A review on the foundations and applications of porous materials can be seen in the pioneering works \cite{Bio72,coussy95,de2000theory,Ehlers2002,Schrefler06,Markert2007}. The fluid volume fraction (porosity) $n_F$ is linked to the fluid volume ratio $\theta$  (fluid content) per unit volume of the undeformed reference configuration $\calB$ via
\begin{equation}
n_F = n_{F,0} + \theta \ ,
\label{porosity-ratio}
\end{equation}
for constant fluid material density, where $n_{F,0}$ is the initial porosity. In the constitutive modeling $\theta$ describes the {\it first local internal} variable (history field). The evolution of this fluid volume ratio $\dot{\theta} = \dot{n}_F$ is derived by the fluid pressure field $p$. The boundary conditions for the pressure are determined as follows
\begin{equation}
p = \bar{p}(\BX,t) \ \textrm{on}\ \partial_D\calB
\AND
\BcalF\cdot\BN = \bar{f}(\BX,t)\ \textrm{on}\ \partial_N\calB\;,
\label{flu-bcs}
\end{equation}
in terms of the material fluid volume flux vector $\BcalF$ and the prescribed fluid pressure $\bar{p}$ and fluid transport $\bar{f}$. The initial condition for the fluid volume ratio is set to $\theta(\BX,t_0)=0$ yields $n_F = n_{F,0}$ in $\calB$. Furthermore, the fluid flux vector in \req{flu-bcs} is linked to the negative material gradient of the fluid pressure via the permeability, according to Darcy-type fluid transport as 
\begin{equation}
\BcalF := - \BK(\BF,d)\; \nabla p \ ,
\end{equation}
where the permeability tensor $\BK$ depends on the material deformation gradient $\BF$ and the crack phase-field $d$. It is decomposed into a {\it Darcy-type flow for an unbroken porous medium} $\boldsymbol{K}_{Darcy}$ and a {\it Poiseuille-type flow in a fully fractured material} $\boldsymbol{K}_{crack}$ defined as
\begin{equation}
\begin{array}{ll}
\boldsymbol{K}(\BF,d) &=
{\boldsymbol{K}_{Darcy}}(\BF)+d^{\zeta}
{\boldsymbol{K}_{frac}}(\BF)
\;, \\ [3mm]
{\boldsymbol{K}_{Darcy}}(\BF) &= \frac{K}{\eta_F} J \BC^{-1}\;, \\ [3mm]
{\boldsymbol{K}_{frac}}(\BF) &= K_c\; \omega^2\; J \big[
\BC^{-1} - \BC^{-1} \Bn \otimes  \BC^{-1} \Bn  
\big]\;,
\end{array}
\end{equation}
as outlined in \cite{aldakheelspp2020,MieheMauthe2015}, where 
$K$ is the intrinsic permeability in an isotropic pore space, $\Bn=\nabla d/|\nabla d|$ is the normal of material crack surface, $\eta_F$ is the dynamic fluid viscosity, $\zeta \ge 1$ is a
permeability transition exponent and $K_c$ is the spatial permeability in fracture. An estimation for the crack width is provided by $\omega = (\lambda_{\bot}-1) h_{e}$ in terms of the stretch perpendicular to the crack $\lambda^2_{\bot} ={\nabla d \cdot \nabla d}/{\nabla d \cdot \BC^{-1} \cdot \nabla d}$ and the characteristic element length $h_{e}$. 

The fluid has to satisfy the balance of fluid mass, reflecting the {\it second} PDE for the coupled problem as
\begin{equation}
\boxed{\dot{n}_F - \bar{r}_F + \Div[\BcalF] = 0 }
\label{2nd-be}
\end{equation}
with a prescribed fluid source $\bar{r}_F$ per unit volume of the reference configuration $\calB$, which describes the injection process in hydraulic fracturing.

For the phase-field problem, a sharp-crack surface topology $\calC \rightarrow \calC_l$ is regularized by the crack surface functional as outlined in \cite{aldakheel+blaz+wriggers18,aldakheel+wriggers+miehe18}
\begin{equation}
\calC_l(d) = \int_{\calB} \gamma_l(d, \nabla d) \, dV
\WITH
\gamma_l(d, \nabla d) =  
\dfrac{1}{2l} d^2 + \dfrac{l}{2} \vert \nabla d \vert^2 \;,
\label{s2-gamma_l}
\end{equation}
based on the crack surface density function $\gamma_l$ per unit volume of the solid and the fracture length scale parameter $l$ that governs the regularization, as plotted in Fig. \ref{Figure1}.
To describe a purely geometric approach to phase-field fracture, the regularized crack phase-field $d$ is obtained by a minimization principle of diffusive crack topology
\begin{equation}
d = \mbox{Arg} \{ \inf_{d}
\calC_l(d) \}
\WITH
d=1 \; \mbox{on} \; \calC\subset \calB\;,
\label{min-d-geo}
\end{equation}
yielding the Euler equation $d - l^2
\Delta d = 0$ in $\calB$ along with the Neumann-type boundary condition $\nabla d \cdot \BN = 0$ on $\partial\calB$. %
Evolution of the regularized crack surface functional \req{s2-gamma_l} can be driven by the constitutive
functions as outlined in \cite{aldakheel16,miehe+hofacker+schaenzel+aldakheel15}, postulating a global evolution equation of regularized crack surface as
\begin{equation}
\vphantom{\int_{\calB}}
\frac{d}{dt} \calC_l(d)  := 
\frac{1}{l} \int_{\calB} [\; (1-d) \calH - \eta  \dot{d}\; ]\;
\dot{d} \, dV \ge 0\;,
\label{gamma-evol}
\end{equation}
where $\eta  \ge 0$ is a material parameter that characterizes 
\textcolor{black}{the artificial/numerical viscosity of the crack propagation}. The crack driving force
\begin{equation}
\calH = \max_{s\in [0,t]} D(\Bx,s) \ge 0
\ ,
\label{driving-force}
\end{equation}
is introduced as the {\it second local history variable} that accounts for the irreversibility of the phase-field evolution by filtering out a maximum value of what is known as the crack driving state function $D$. Then the evolution statement
\req{gamma-evol} provides the local equation for the
evolution of the crack phase-field in the domain $\calB$ along with its homogeneous Neumann boundary condition as
\begin{equation}
\boxed{
[ \, d - l^2 \Delta d \, ] + \eta \dot{d} + (d-1) {\calH} = 0
}
\label{euler-eq-d}
\end{equation}
with $\nabla d \cdot \BN = 0$ on $\partial\calB$. It represents the {\it third} PDE for the coupled problem.

\sectpb{Constitutive functions}

The multi-physics problem is based on three primary fields to characterize the hydro-poro-elasticity of fluid-saturated porous media as
\begin{equation}
\mbox{Global Primary Fields}: \BfrakU := \{ \Bvarphi, p, d \}
\label{global-fields}
\ ,
\end{equation}
the deformation map $\Bvarphi$, the pressure field $p$ and the crack phase-field $d$. The constitutive approach to hydraulic phase-field fracture in poroelastic media focuses on the set
\begin{equation}
\mbox{Constitutive State Variables}: 
\BfrakC := \{ \BF, \theta, d, \nabla d \}
\ ,
\label{state}
\end{equation}
reflecting a combination of poro-elasticity with a first-order gradient damage modeling. It is based on the definition of a pseudo-energy density per unit volume contains the sum
\begin{equation}
\boxed{{W}(\BfrakC) = 
	{W}_{elas}(\BF, d) + {W}_{fluid}(\BF,\theta) + 
	{W}_{frac}(d, \nabla d)
	}
\label{pseudo-energy}
\end{equation}
of a degrading elastic part ${W}_{elas}$ and a contribution due to fluid ${W}_{fluid}$ and fracture ${W}_{frac}$ that contain the accumulated dissipative energy. The elastic contribution is modeled with a Neo-Hookean strain energy function for a homogeneous compressible isotropic elastic solid
\begin{equation}
{W}_{elas}(\BF, d)= g(d)\; {\psi}_{elas}(\BF) 
\WITH
{\psi}_{elas}(\BF) = \frac{\mu}{2}\Big[(\BF:\BF -3) + \frac{2}{\beta} (J^{-\beta} -1) \Big]\;,
\label{elas-part}
\end{equation} 
in terms of the shear modulus $\mu$ and the parameter $\beta = 2\nu/(1-2\nu)$ with the Poisson number $\nu$. The function ${g}(d) = (1-d)^2$ models the degradation of the elastic energy of the solid due to fracture. It interpolates between the unbroken response for $d=0$ and the fully broken state at $d=1$ by satisfying the constraints
${g}(\mathit{0}) = 1$,
${g}(\mathit{1}) = 0$,
${g}{}^{\prime}(d) \le 0$ and
${g}{}^{\prime}(\mathit{1}) = 0$.
The fluid contribution is assumed to have the form
\begin{equation}
{W}_{fluid}(\BF, \theta) = \frac{M}{2} \Bigg[ B^2 (J-1)^2 - 2\, \theta\,B (J-1) +  \theta^2
\Bigg]\;,
\label{fluid-part}
\end{equation}
in terms of the Biot's coefficient $B$ and Biot's modulus $M$. Following the Coleman-Noll procedure, the fluid pressure $p$ and the first Piola-Kirchoff stress tensor $\BP$ are obtained from the pseudo-energy density function ${W}$ in \req{pseudo-energy} for isotropic material behavior as
\begin{equation}
\begin{array}{ll}
p(\BF, \theta) &:= \frac{\partial{W}}{\partial \theta } = \theta M  - M B (J-1)\;,
\\[4mm]
\BP(\BF, p, d) &:= \frac{\partial{W}}{\partial\BF }= g(d) \BP_{eff}(\BF) -
B p J \BF^{-T}
\WITH
\BP_{eff} = \mu \big[ \BF - J^{-\beta} \BF^{-T}
\big]\;,
\end{array}
\label{p-piola-stresses}
\end{equation}
where the stress tensor is additively decomposed into an effective part $\BP_{eff}$ and a pressure part according to the classical Terzaghi split, as outlined in \cite{terzaghi1943theoretical,de1990development}. Using the pressure definition in \reqi{p-piola-stresses}1 and the second PDE in \req{2nd-be} along with \req{porosity-ratio}, the balance of mass is modified as follows
\begin{equation}
\frac{\dot{p}}{M} + B \dot{J} - \bar{r}_F+ \Div[\BcalF] = 0\;,
\label{pres-pde}
\end{equation}
which now depends on the fluid pressure $p$ and the deformation $\Bvarphi$. 

The fracture part of pseudo-energy density \req{pseudo-energy} is modeled by
\begin{equation}
{W}_{frac}(d, \nabla d) =  [1 - {g}(d)]\; \psi_c +
2 {\psi_c}\; l \;{\gamma}_l(d, \nabla d)\;,
\label{frac-part}
\end{equation}
where ${\psi}_c > 0$ is a critical fracture energy. It is defined in terms of the critical effective stress $\sigma_c$ or the Griffith's energy release rate $G_c$, as outlined in \cite{aldakheel+blaz+wriggers18}
\begin{equation}
\psi_c = \frac{\sigma_c^2}{2E} = \frac{3}{8 l \sqrt{2}} G_c \; .
\end{equation}
By taking the variational derivative $\delta_d W$ of \req{pseudo-energy} with some manipulation as documented in \cite{aldakheeletal18}, the third PDE in \req{euler-eq-d} yields for the rate-independent setting as follows
\begin{equation}
2\psi_c [d - l^2 \Delta d] + 2 (d-1) \calH = 0\;, 
\label{frac-eqs}
\end{equation}
in terms of the history field $\calH$, introduced in \req{driving-force}. The crack driving state function $D$ is defined by
\begin{equation}
D:= \Big< \psi_{elas}(\BF(\BX,s)) - \psi_c \Big>_+ \ge 0\;,
\label{h-history-field}
\end{equation}
with the Macaulay bracket $\langle x \rangle_+ := (x + \vert x\vert)/2$, that ensures the irreversibility of the crack evolution.

\sectpb{Weak formulations for the coupled problem}
The update of the primary fields $\BfrakU$ in \req{global-fields} in a typical time increment $[t_n, t_{n+1}]$ with time step $\Delta t >0$ is governed by three PDEs in \req{equil:defo}, \req{pres-pde} and \req{euler-eq-d} in a strong form setting. Next, we define three test functions for the deformation $\delta \Bvarphi(\BX) \in\{ {\BH}^1(\calB)^\delta: \delta\Bvarphi=\Bzero \; \mathrm{on} \; \partial_D\calB \}$, fluid pressure $\delta p(\BX) \in\{ {H}^1(\calB): \delta p=0 \; \mathrm{on} \; \partial_D\calB \}$ and crack phase-field $\delta d(\BX)\in {H}^1(\calB)$. The weak formulations for the above introduced three PDEs of the coupled poro-elastic media problem at fracture are derived from a standard Galerkin procedure as
\begin{equation}
\begin{array}{ll}
G_\varphi(\BfrakU, \delta \Bvarphi) &= \int_\calB \Big[ \BP:\nabla \delta \Bvarphi - \bar{\Bb} \cdot \delta \Bvarphi \Big] dV - \int_{\partial_N\calB} \bar{\Btau} \cdot \delta \Bvarphi  \; dA
= 0 \ , \\ [4mm]
G_p(\BfrakU, \delta p) &= \int_\calB \Big[\Big(\frac{1}{M}(p-p_n) + B (J-J_n) -\Delta t \;\bar{r}_F
\Big)\delta p + (\Delta t \;\BK \;\nabla p) \cdot \nabla \delta p \Big] dV \\ [3mm]
& + \int_{\partial_N\calB} \bar{f} \;\delta p\; dA = 0 \ , \\ [4mm]
G_d(\BfrakU, \delta d) &= \int_\calB \Big[ \Big(2\psi_c \;d + 2(d-1) \calH \Big) \delta d  + 2\psi_c \;l^2 \; \nabla d \cdot \nabla \delta d
\Big] dV = 0 \ .
\end{array}
\label{weakForm}
\end{equation}
This set of equations describes the constitutive model fully. Next, we use \req{weakForm} as a departure point for the global-local approach in Section \ref{Section3}.
\sectpa[Section3]{Extension Towards Global-Local Formulations}

\begin{figure}[!t]
	\centering
	{\includegraphics[clip,trim=8cm 6cm 5cm 13cm, width=16cm]{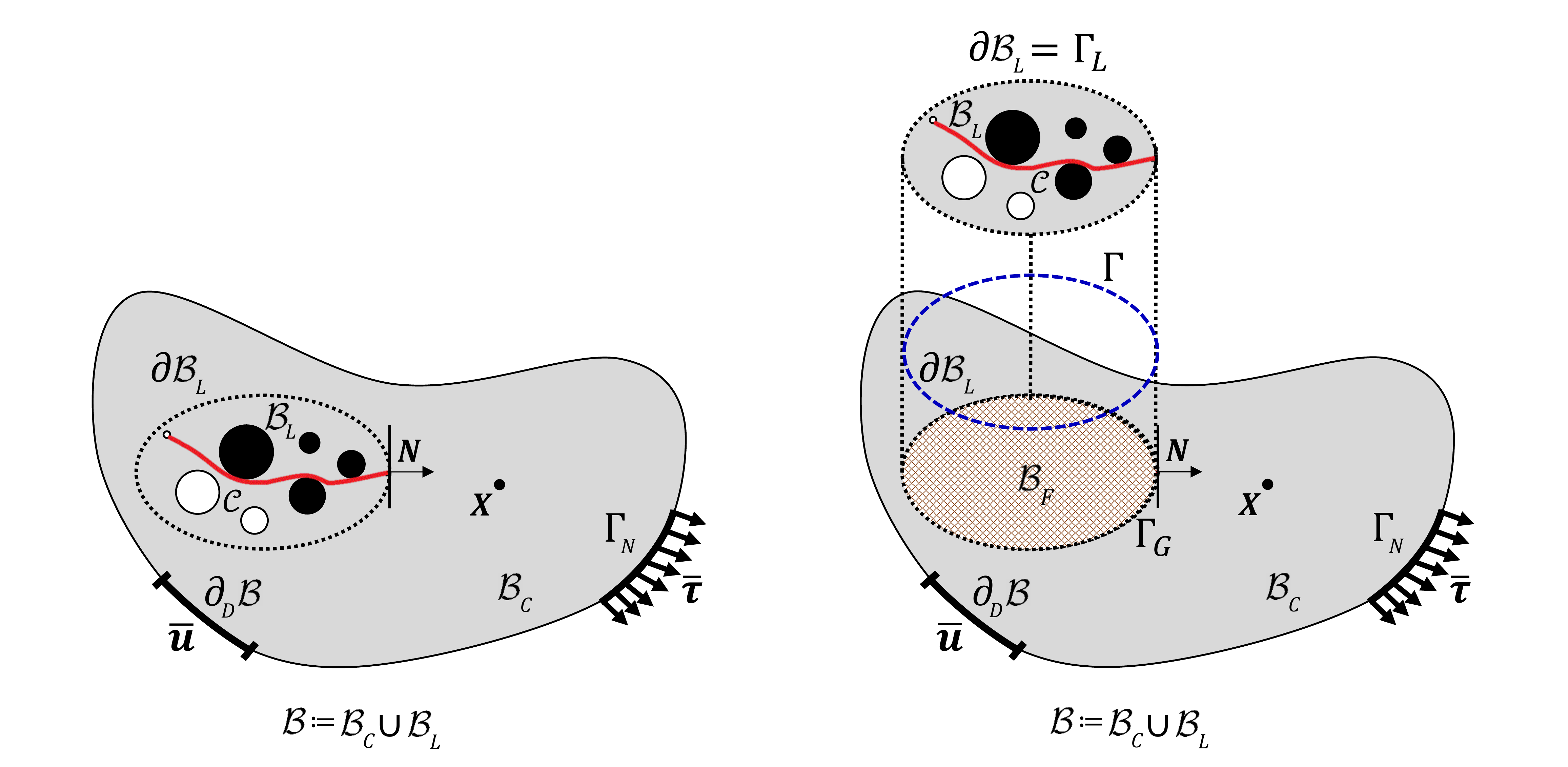}}  
	\caption{Illustration of the Global-Local formulation. $(a)$ Geometry and loading setup of the single-scale boundary value problem. $(b)$ Global-Local setting, by introduction of the fictitious domain $\calB_F$ through prolongation of $\calB_C$ to the entire domain whereas its unification is so-called global domain $\calB_G:=\calB_C\cup\Gamma\cup\calB_F$.
	}
	\label{Fig3}
\end{figure}

In this section the above introduced system of equations for the coupled problem will be solved using the Global-Local (GL) method, that is rooted in the domain decomposition approach \cite{DDRey06}. It represents an initial contribution to the use of the GL formulation at {\it large deformations} for solving fracture mechanics problems numerically. The main objective here is to introduce an adoption of the hydraulic phase-field fracture formulation in poroelastic media within legacy codes, specifically for industrial applications. To this end, the material body $\calB$ is decomposed into a global domain $\calB_G$ representing a poro-elastic media and a local domain $\calB_L$ reflecting the hydraulic fracturing (fracking) region. The global domain $\calB_G:= \calB_C\cup\calB_f\cup\Gamma$ is further split into a complementary domain $\calB_C$ corresponds to the intact area, a fictitious domain $\calB_f$ depicts a coarse projection of the local domain into the global one and an interface $\Gamma$
between the unfractured and the fractured domains. The {fictitious domain} $\calB_{f}$ is a prolongation of $\calB_{C}$ towards ${\calB}$, i.e. {\it recovering the space} of ${\calB}$ that is obtained by removing $\calB_{L}$ from its continuum domain, see Fig. \ref{Fig3}. This gives the same constitutive modeling used in $\calB_{C}$ for $\calB_{f}$. We also use the identical discretization space for both $\calB_{f}$ and $\calB_{C}$, which results in $h_f:=h_C$. The external loads are applied on $\calB_{C}$ and hence $\calB_{L}$ is assumed to be free from external loads. Such assumption is standard for the multi-scale setting, see \cite{Fish2014}.

At the interface $\Gamma$, global and local interfaces denoted as ${\Gamma}_{G}\subset\calB_{G}$ and $\Gamma_{L}\subset\calB_{L}$ are defined, such that in the continuum setting we have ${\Gamma}=\Gamma_{G}=\Gamma_{L}$. Hence, the deformation map $\Bvarphi$ and the fluid pressure $p$ for both global and local domains do exactly coincide in the strong sense at interface, yielding
\begin{equation}
{\Bvarphi}_L(\BX,t)\overset{!}{=}{\Bvarphi}_G(\BX,t) 
\AND
p_L(\BX,t)\overset{!}{=} p_G(\BX,t) 
\quad  \mbox{at} \quad {\BX}\in\Gamma\;.
\label{uCont}
\end{equation}
However in a discrete setting we might have ${\Gamma}\neq\Gamma_{G}\neq\Gamma_{L}$ due to the presence of different meshing schemes (i.e. different element size/type used in $\calB_{G}$ and $\calB_{L}$ such that $h \neq h_L \neq h_G$ on ${\Gamma}$). 
\begin{Remark}
	\label{non_match}
The strong deformation/pressure continuity requirement given in Eq. \ref{uCont} is too restrictive from the computational standpoint \cite{Farhat1991}. To resolve the phase field problem, one requires {$h_L \ll h_G$}. However, if we assume ${\Bvarphi}_L\overset{!}{=}{\Bvarphi}_G$ and $p_L\overset{!}{=}p_G$ on $\Gamma$, this yields $\Gamma_{L}=\Gamma_{G}$ in a discretized setting hence $h_L = h_G$ on $\Gamma$ which contradicts {$h_L \ll h_G$}.
\end{Remark}
\begin{form}[Continuity conditions at interface]
	\label{form_interface}
	Following Remark \ref{non_match}, we modify Eq. \ref{uCont} in a strong sense by introducing the deformation $\Bvarphi_\Gamma(\BX,t)$ and pressure $p_\Gamma(\BX,t)$ interface and their corresponding traction forces $\{\Blambda_L^\varphi, \Blambda^\varphi_C\}$ and $\{\lambda_L^p, \lambda_C^p\}$ that are introduced as Lagrange multipliers. This results in a set of equations at the interface as
	\begin{align*}
	\left\{
	\begin{tabular}{ll}
	${\Bvarphi}_L(\BX,t)  = {\Bvarphi}_{\Gamma}(\BX,t) $ & \mbox{at} ${\BX}\in{\Gamma}_{L}$, \\[0.1cm]
	${\Bvarphi}_G(\BX,t)  = {\Bvarphi}_{\Gamma}(\BX,t) $ & \mbox{at} ${\BX}\in{\Gamma}_{G}$, \\[0.1cm]
	$\Blambda_L^\varphi(\BX,t) + \Blambda^\varphi_C(\BX,t) = \Bzero$ & \mbox{at} ${\BX}\in{\Gamma}$,
	\end{tabular}
	\right.
	\ \mbox{and} \
	\left\{
	\begin{tabular}{ll}
	$p_L(\BX,t) = p_{\Gamma}(\BX,t)$ & \mbox{at} ${\BX}\in{\Gamma}_{L}$, \\[0.1cm]
	$p_G(\BX,t) = p_{\Gamma}(\BX,t)$ & \mbox{at} ${\BX}\in{\Gamma}_{G}$, \\[0.1cm]
	$\lambda_L^p(\BX,t) + \lambda^p_C(\BX,t) = \mathit{0}$ & \mbox{at} ${\BX}\in{\Gamma}$.
	\end{tabular}
	\right.
	\end{align*}
\end{form}
\bigskip
Accordingly, the {\it single-scale} deformation map ${\Bvarphi}(\BX,t)$ and fluid pressure field $p(\BX,t)$ in Section \ref{Section2} are decomposed as follows
\begin{equation}
{\Bvarphi}(\BX,t) =\left\{
\begin{tabular}{ll}
${\Bvarphi}_L(\BX,t)$ & for ${\BX}\in\calB_{L}$, \\[0.1cm]
${\Bvarphi}_G(\BX,t)$ & for ${\BX}\in\calB_{G}$, \\[0.1cm]
${\Bvarphi}_{\Gamma}(\BX,t)$ & for ${\BX}\in{\Gamma}$,
\end{tabular}
\right.
\AND
p(\BX,t) =\left\{
\begin{tabular}{ll}
$p_L(\BX,t)$ & for ${\BX}\in\calB_{L}$, \\[0.1cm]
$p_G(\BX,t)$ & for ${\BX}\in\calB_{G}$, \\[0.1cm]
$p_{\Gamma}(\BX,t)$ & for ${\BX}\in{\Gamma}$.
\end{tabular}
\right.
\label{phi-p-DD}
\end{equation}
The fracture surface lives only in $\calB_{L}$. Hence we can introduce scalar-valued function $d_L(\BX,t):\calB_{L}\rightarrow[0,1]$. The {\it single-scale phase-field} $d$ is then decomposed in the following representation
\begin{equation}
d(\BX,t)  :=\left\{
\begin{tabular}{ll}
$d_L$ & for ${\BX}\in\calB_{L}$, \\[0.1cm]
$0$ & for ${\BX}\in\calB_{G}$.
\end{tabular}
\right.
\label{dL-DD}
\end{equation}
Now the multi-physics problem for the Global-Local approach is based on eleven primary fields to characterize the hydro-poro-elasticity of fluid-saturated porous media at finite strains as
\begin{equation}
\mbox{Extended Primary Fields}: \BfrakP := \{ \Bvarphi_G, \Bvarphi_L, p_G, p_L, d_L, \Blambda^\varphi_C, \Blambda^\varphi_L, \lambda^p_C, \lambda^p_L, \Bvarphi_\Gamma, p_\Gamma \}
\label{gl-fields}
\ .
\end{equation}
%
\sectpb[Section31]{Governing formulations for the Global-Local coupling system}
Based on the above introduced decompositions and the weak formulations outlined in \req{weakForm}, this section describes the GL weak forms of the PDEs for the coupled problem. The global weak formulations of the deformation and pressure field take the form
\begin{align*}\label{GD}
G_{\Bvarphi_G}(\BfrakP, \delta\Bvarphi_G)&:=
\int_{\calB_{G}} \BP(\nabla \Bvarphi_G, p_G, 0):\nabla \delta\Bvarphi_G \mathrm{d}V - \int_{\calB_{f}} \BP(\nabla \Bvarphi_G, p_G, 0):\nabla \delta \Bvarphi_G \mathrm{d}V\\
&-\int_{\Gamma_{G}} \Blambda_C^\varphi \cdot \delta {{\Bvarphi}}_G\,\mathrm{d}A
-\int_{\Gamma_{N,G} } {{\bar\Btau}} \cdot \delta{{\Bvarphi}}_G\,\mathrm{d}A=0 \ ,\\[1mm]
G_{p_G}(\BfrakP, \delta p_G)&:= \int_{\calB_G} \Big[\frac{1}{M}(p_G-p_{G,n}) + B \Big(J(\nabla \Bvarphi_G)-J_n(\nabla \Bvarphi_G)\Big)  \Big]\delta p_G\; \mathrm{d}V\\
\tag{G}&
+ \int_{\calB_G} \Big[ (\Delta t \;\BK(\nabla \Bvarphi_G, 0) \;\nabla p_G) \cdot \nabla \delta p_G \Big] \mathrm{d}V\\
&+ \int_{\calB_f} \Big[\frac{1}{M}(p_G-p_{G,n}) + B \Big(J(\nabla \Bvarphi_G)-J_n(\nabla \Bvarphi_G)\Big) -\Delta t \;\bar{r}_F
\Big]\delta p_G\; \mathrm{d}V\\
&
+ \int_{\calB_f} \Big[ (\Delta t \;\BK(\nabla \Bvarphi_G, 0) \;\nabla p_G) \cdot \nabla \delta p_G \Big] \mathrm{d}V\\
& -\int_{\Gamma_{G}} \lambda_C^p \ \delta {p}_G\,\mathrm{d}A
+ \int_{\Gamma_{N,G}} \bar{f} \;\delta p_G\; dA = 0 \ ,
\end{align*}
where $\delta \Bvarphi_G \in\{ {\BH}^1(\calB_G)^\delta: \delta\Bvarphi_G=\Bzero \; \mathrm{on} \; \partial_D\calB \}$ and $\delta p_G \in\{ {H}^1(\calB_G): \delta p_G=0 \; \mathrm{on} \; \partial_D\calB \}$ are the global test functions. Note that the pressure injection process of hydraulic fracturing $\bar{r}_F$ {\it exists only} in the fictitious domain $\calB_f$. The local weak formulations assumes the form
\begin{align*}\label{LD}
G_{\Bvarphi_L}(\BfrakP, \delta\Bvarphi_L)&:=
\int_{\calB_{L}} \BP(\nabla \Bvarphi_L, p_L, d_L):\nabla \delta\Bvarphi_L \mathrm{d}V - \int_{\Gamma_{L}} \Blambda_L^\varphi \cdot \delta {{\Bvarphi}}_L\,\mathrm{d}A=0 \ ,\\[1mm]
\tag{L}
G_{p_L}(\BfrakP, \delta p_L)&:= \int_{\calB_L} \Big[\frac{1}{M}(p_L-p_{L,n}) + B \Big(J(\nabla \Bvarphi_L)-J_n(\nabla \Bvarphi_L)\Big)  \Big]\delta p_L\; \mathrm{d}V\\
&
+ \int_{\calB_L} \Big[ \Big(\Delta t \;\BK(\nabla \Bvarphi_L, d_L) \;\nabla p_L \Big) \cdot \nabla \delta p_L \Big] \mathrm{d}V
-\int_{\Gamma_{L}} \lambda_L^p \ \delta {p}_L\,\mathrm{d}A = 0\\
G_{d_L}(\BfrakP, \delta d_L)&: = \int_{\calB_L} \Big[ \Big(2\psi_c \;d_L + 2(d_L-1) \ \calH(\nabla \Bvarphi_L) \Big) \delta d_L  + 2\psi_c \;l^2 \; \nabla d_L \cdot \nabla \delta d_L
\Big] dV = 0 \ ,
\end{align*}
where $\delta \Bvarphi_L \in {\BH}^1(\calB_L)$, $\delta p_L \in {H}^1(\calB_L)$ and $\delta d_L \in {H}^1(\calB_L)$ are the local test functions for the deformation, fluid pressure and crack phase-field, respectively. 

Next, we derive the weak formulations for the deformation and pressure continuity at interface ${\Gamma}$ introduced in Formulation \ref{form_interface} by using a standard Galerkin procedure
\begin{equation*}
G_{\Bvarphi_\Gamma}(\BfrakP, \delta \Bvarphi_\Gamma):=\int_\Gamma (\Blambda^\varphi_C+\Blambda^\varphi_L) \cdot \delta \Bvarphi_\Gamma \,\mathrm{d}A=0,
\label{Coupl1}
\tag{C$_1$}
\end{equation*}
\begin{equation*}
G_{\Blambda^\varphi_C}(\BfrakP, \delta \Blambda^\varphi_C):=\int_\Gamma (\Bvarphi_\Gamma-\Bvarphi_G) \cdot \delta \Blambda^\varphi_C \,\mathrm{d}A=0,
\label{Coupl2}
\tag{C$_2$}
\end{equation*}
\begin{equation*}
G_{\Blambda^\varphi_L}(\BfrakP, \delta \Blambda^\varphi_L):=\int_\Gamma (\Bvarphi_\Gamma-\Bvarphi_L) \cdot \delta \Blambda^\varphi_L \,\mathrm{d}A=0,
\label{Coupl3}
\tag{C$_3$}
\end{equation*}
\begin{equation*}
G_{p_\Gamma}(\BfrakP, \delta p_\Gamma):=\int_\Gamma (\lambda^p_C+\lambda^p_L)  \delta p_\Gamma \,\mathrm{d}A=0,
\label{Coupl4}
\tag{C$_4$}
\end{equation*}
\begin{equation*}
G_{\lambda^p_C}(\BfrakP, \delta \lambda^p_C):=\int_\Gamma (p_\Gamma-p_G) \delta \lambda^p_C \,\mathrm{d}A=0,
\label{Coupl5}
\tag{C$_5$}
\end{equation*}
\begin{equation*}
G_{\lambda^p_L}(\BfrakP, \delta \lambda^p_L):=\int_\Gamma (p_\Gamma-p_L) \delta \lambda^p_L \,\mathrm{d}A=0,
\label{Coupl6}
\tag{C$_6$}
\end{equation*}
herein $\delta \Bvarphi_\Gamma\in{\BH}^1(\Gamma)$; $\delta p_\Gamma\in{H}^1(\Gamma)$; $\delta \Blambda^\varphi_C,\delta\Blambda^\varphi_L \in{\BL}^2(\Gamma)$ and $ \delta\lambda^p_C, \delta\lambda^p_L \in{L}^2(\Gamma)$ are the corresponding test functions. Equations (G), (L) and (C1)--(C6) specify the entire system of the Global-Local approach.

\sectpb[Section32]{Dirichlet-Neumann type boundary conditions}
Within a Global-Local computational scheme, instead of finding the stationary solution of the (\ref{GD}), (\ref{LD}) along with (\ref{Coupl1}) -- (\ref{Coupl6}) in the monolithic sense, an alternate minimization is employed. {This} is in line with \cite{NoiiGL18}, which leads to the Global-Local formulation using the concept of non-intrusiveness. Here the global and local level are solved in a multiplicative manner according to the idea of Schwarz' alternating method \cite{alternatMota17}.

Let $k\geq0$ be the Global-Local iteration index at a fixed loading step $n$. The iterative solution procedure for the Global-Local computational scheme is as {follows:}
\begin{itemize}
	\item Dirichlet local problem: solution of local problem (\ref{LD}) coupled with (\ref{Coupl3}) and (\ref{Coupl6}),
	\item Pre-processing global level: recovery phase using (\ref{Coupl1}) and (\ref{Coupl4}),
	\item Neumann global problem: solution of global problem (\ref{GD}),
	\item Post-processing global level: recovery phase using (\ref{Coupl2}) and (\ref{Coupl5}).
\end{itemize}
Despite of its strong non-intrusiveness \cite{Allix09}, there are two shortcomings embedded in the system which have to be resolved. ({a}) Due to the extreme difference in stiffness between the local domain and its projection to the global level, i.e. fictitious domain $\calB_f$, a relaxation/acceleration techniques has to be used, see \cite{NoiiGL18}. ({b})~Additionally, it turns out that if the solution vector  
\begin{equation}
\BfrakP^k = (\Bvarphi_G^k, \Bvarphi_L^k, p_G^k, p_L^k, d_L^k, \Blambda^{\varphi,k}_C, \Blambda^{\varphi,k}_L, \lambda^{p,k}_C, \lambda^{p,k}_L, \Bvarphi_\Gamma^k, p_\Gamma^k)
\end{equation}
is plugged into equations (\ref{GD}), (\ref{LD}), (\ref{Coupl1}) -- (\ref{Coupl6}), the imbalanced quantities for the deformation and pressure fields follow
\begin{equation}
\int_\Gamma (\Bvarphi_\Gamma^k-\Bvarphi_L^k) \cdot \delta {\Blambda}_L^\varphi  \,\mathrm{d}{A} \neq 0 
\AND
\int_\Gamma (p_\Gamma^k-p_L^k) \; \delta {\lambda}_L^p  \,\mathrm{d}{A} \neq 0 ,
\label{t2x}
\end{equation}
resulting in the {\it iterative} Global-Local computation scheme. The aforementioned difficulties motivate us to provide alternative coupling conditions that overcome these challenges, which are explained in the following section.

\sectpb[Section33]{Robin-type boundary conditions}
In this section, the Global-Local formulation is enhanced using Robin-type boundary conditions to relax the stiff local response that is observed at the global level (due to the local non-linearity). Furthermore the computational time is reduced. This improves the resolution of the imbalanced quantities in (\ref{t2x}) and it accelerates Global-Local computational iterations.

Recall, the coupling equations denoted in (\ref{Coupl1}) -- (\ref{Coupl6}) arise from the continuity conditions at the interface in a strong sense. That provides the boundary conditions which have to be imposed on the global and local levels. At that level the Robin-type boundary conditions are formulated.

\sectpc[Section331]{Robin-type boundary conditions at the local level}

\sectpd[Section3311]{Finite deformation}
For the mechanical deformation field at the local level, a new coupling term is introduced as a combination of (\ref{Coupl1}) and (\ref{Coupl2})
\begin{equation}
G_{\Bvarphi_\Gamma}(\BfrakP, \delta \Bvarphi_\Gamma) + \KIA_L^\varphi
G_{\Blambda^\varphi_C}(\BfrakP, \delta \Blambda^\varphi_C) = 
\int_\Gamma (\Blambda^\varphi_C+\Blambda^\varphi_L) \cdot \delta \Bvarphi_\Gamma \,\mathrm{d}A + 
\KIA_L^\varphi \int_\Gamma (\Bvarphi_\Gamma-\Bvarphi_G) \cdot \delta \Blambda^\varphi_C \,\mathrm{d}A = 0
\label{RBC1}
\end{equation}
This leads for iteration $k$ to
\begin{equation}
\int_\Gamma (\Blambda^{\varphi,k-1}_C+\Blambda^{\varphi,k}_L) \cdot \delta {\Bvarphi}_\Gamma \,\mathrm{d}A+\KIA_L^{\varphi}\int_\Gamma (\Bvarphi_\Gamma^{k,\frac{1}{2}}-\Bvarphi^{k-1}_G) \cdot \delta {\Blambda}^{\varphi}_C \,\mathrm{d}A=0.
\label{RBC2}
\end{equation}
Herein, $\KIA_L^\varphi$ is a local augmented stiffness matrix for the deformation applied at the interface which serves as regularization of the local Jacobian matrix. By means of (\ref{RBC2}) at iteration $k$, the local system of equations for the mechanical problem at the interface (\ref{Coupl1}) -- (\ref{Coupl3}) results in the following modified boundary conditions
\begin{equation*}
\int_\Gamma \Blambda^{\varphi,k}_L\cdot \delta {\Bvarphi}_\Gamma \,\mathrm{d}A+\KIA_L^\varphi\int_\Gamma \Bvarphi_\Gamma^{k,\frac{1}{2}} \cdot \delta {\Blambda}_C^\varphi \,\mathrm{d}A={\BLambda}^{\varphi,k-1}_L,
\label{RBC3}
\tag{$\widetilde {\text{C}}_1$}
\end{equation*}
\begin{equation*}
\int_\Gamma (\Bvarphi_\Gamma^{k,\frac{1}{2}}-\Bvarphi^{k}_L) \cdot \delta {\Blambda}_L^\varphi \,\mathrm{d}A=0,
\label{RBC4}
\tag{$\widetilde {\text{C}}_2$}
\end{equation*}  
with
\begin{equation}
{\BLambda}^{\varphi,k-1}_L:=\BLambda_L(\Blambda^{\varphi,k-1}_C,\Bvarphi_G^{k-1};\KIA_L^\varphi)=
\KIA_L^\varphi\int_\Gamma \Bvarphi_G^{k-1} \cdot \delta {\Blambda}_C^\varphi \,\mathrm{d}A-\int_\Gamma \Blambda^{\varphi,k-1}_C\cdot \delta {\Bvarphi}_\Gamma \,\mathrm{d}A \ .
\label{RBC5}
\end{equation}

\sectpd[Section3311]{Fluid pressure}
Analogously to the coupling terms for the deformation introduced above, we modify the local system of equations for the pressure field at the interface (\ref{Coupl4}) -- (\ref{Coupl6}). It results in the following modified boundary conditions
\begin{equation*}
\int_\Gamma \lambda^{p,k}_L \;\delta {p}_\Gamma \,\mathrm{d}A+\KIA_L^p\int_\Gamma p_\Gamma^{k,\frac{1}{2}}\; \delta {\lambda}_C^p \,\mathrm{d}A={\Lambda}^{p,k-1}_L,
\label{RBC6}
\tag{$\widetilde {\text{C}}_3$}
\end{equation*}
\begin{equation*}
\int_\Gamma (p_\Gamma^{k,\frac{1}{2}}-p^{k}_L) \delta {\lambda}_L^p \,\mathrm{d}A=0,
\label{RBC7}
\tag{$\widetilde {\text{C}}_4$}
\end{equation*}  
with
\begin{equation}
{\Lambda}^{p,k-1}_L:=\Lambda_L(\lambda^{p,k-1}_C,p_G^{k-1};\KIA_L^p)=
\KIA_L^p\int_\Gamma p_G^{k-1} \delta {\lambda}_C^p \,\mathrm{d}A-\int_\Gamma \lambda^{p,k-1}_C \delta {p}_\Gamma \,\mathrm{d}A.
\label{RBC8}
\end{equation}
Along with (\ref{LD}), the local system of equations has to be solved for $(\Bvarphi^{k}_L, p^{k}_L, \Blambda^{\varphi,k}_L, \lambda^{p,k}_L, \Bvarphi_\Gamma^{k,\frac{1}{2}}, p_\Gamma^{k,\frac{1}{2}})$ for given local Robin-type parameters $({\BLambda}^{\varphi,k-1}_L, {\Lambda}^{p,k-1}_L, \KIA_L^\varphi, \KIA_L^p)$.

\begin{Remark}
	\label{k-iteration-choice}
	In the numerical implementation, the current local fields are computed based on the old global variables as history fields, see \req{RBC2}. Hereby, the deformation $\Bvarphi_\Gamma$ and fluid pressure $p_\Gamma$ at the interface are updated at iteration ($k,\frac{1}{2}$). This choice is essential for the construction of the Robin-type boundary conditions. Note that, we proved in previous work that $\Bu_\Gamma^{(k,\frac{1}{2})} = \Bu_\Gamma^k$ with $\Bvarphi:= \Bu + \BX$ where $\BX$ is a fixed initial configuration, see \cite{NoiiAldakheelWickWriggers2019}. With this prove at interface the continuity conditions are satisfied yielding a well posed problem and accelerate the convergence results. Note that other coupling conditions at the interface, i.e. updating the deformation and pressure at iteration $k$ in ({$\widetilde {\text{C}}_2$}) and ({$\widetilde {\text{C}}_4$})  gives ill-posed problem due to the imposition of both Neumann and Dirichlet boundary conditions at same time at $\Gamma$.  
\end{Remark}

\sectpc[Section332]{Robin-type boundary conditions at the global level}

\sectpd[Section3321]{Finite deformation}
Accordingly, at the global level, the new coupling term is stated as a combination of (\ref{Coupl1}) and (\ref{Coupl3}) for the mechanical deformation as
\begin{equation}
G_{\Bvarphi_\Gamma}(\BfrakP, \delta \Bvarphi_\Gamma) + \KIA_G^\varphi
G_{\Blambda^\varphi_L}(\BfrakP, \delta \Blambda^\varphi_L) = 
\int_\Gamma (\Blambda^\varphi_C+\Blambda^\varphi_L) \cdot \delta \Bvarphi_\Gamma \,\mathrm{d}A + 
\KIA_G^\varphi \int_\Gamma (\Bvarphi_\Gamma-\Bvarphi_L) \cdot \delta \Blambda^\varphi_L \,\mathrm{d}A = 0
\label{RBC9}
\end{equation}
This leads for iteration $k$ to
\begin{equation}
\int_\Gamma (\Blambda^{\varphi,k}_C+\Blambda^{\varphi,k}_L) \cdot \delta {\Bvarphi}_\Gamma \,\mathrm{d}A+\KIA_G^{\varphi}\int_\Gamma (\Bvarphi_\Gamma^{k}-\Bvarphi^{k}_L) \cdot \delta {\Blambda}^{\varphi}_L \,\mathrm{d}A=0 \ ,
\label{RBC10}
\end{equation}
where, $\KIA_G^{\varphi}$ is a global augmented stiffness matrix for the deformation applied on the interface. Through (\ref{RBC10}) at the iteration $k$, the Robin-type boundary condition at the global level follows
\begin{equation*}
\int_\Gamma \Blambda^{\varphi,k}_C\cdot \delta {\Bvarphi}_\Gamma \,\mathrm{d}A+\KIA_G^\varphi\int_\Gamma \Bvarphi_\Gamma^{k} \cdot \delta {\Blambda}_L^\varphi \,\mathrm{d}A={\BLambda}^{\varphi,k}_G,
\label{RBC11}
\tag{$\widetilde {\text{C}}_5$}
\end{equation*}
\begin{equation*}
\int_\Gamma (\Bvarphi_\Gamma^{k,\frac{1}{2}}-\Bvarphi^{k}_G) \cdot \delta {\Blambda}_C^\varphi \,\mathrm{d}A=0,
\label{RBC12}
\tag{$\widetilde {\text{C}}_6$}
\end{equation*}  
with
\begin{equation}
{\BLambda}^{\varphi,k}_G:=\BLambda_G(\Blambda^{\varphi,k}_L,\Bvarphi_L^{k};\KIA_G^\varphi)=
\KIA_G^\varphi\int_\Gamma \Bvarphi_L^{k} \cdot \delta {\Blambda}_L^\varphi \,\mathrm{d}A-\int_\Gamma \Blambda^{\varphi,k}_L\cdot \delta {\Bvarphi}_\Gamma \,\mathrm{d}A \ .
\label{RBC13}
\end{equation}

\begin{figure}[!b]
	\centering
	{\includegraphics[clip,trim=0cm 3cm 0cm 4.5cm, width=16cm]{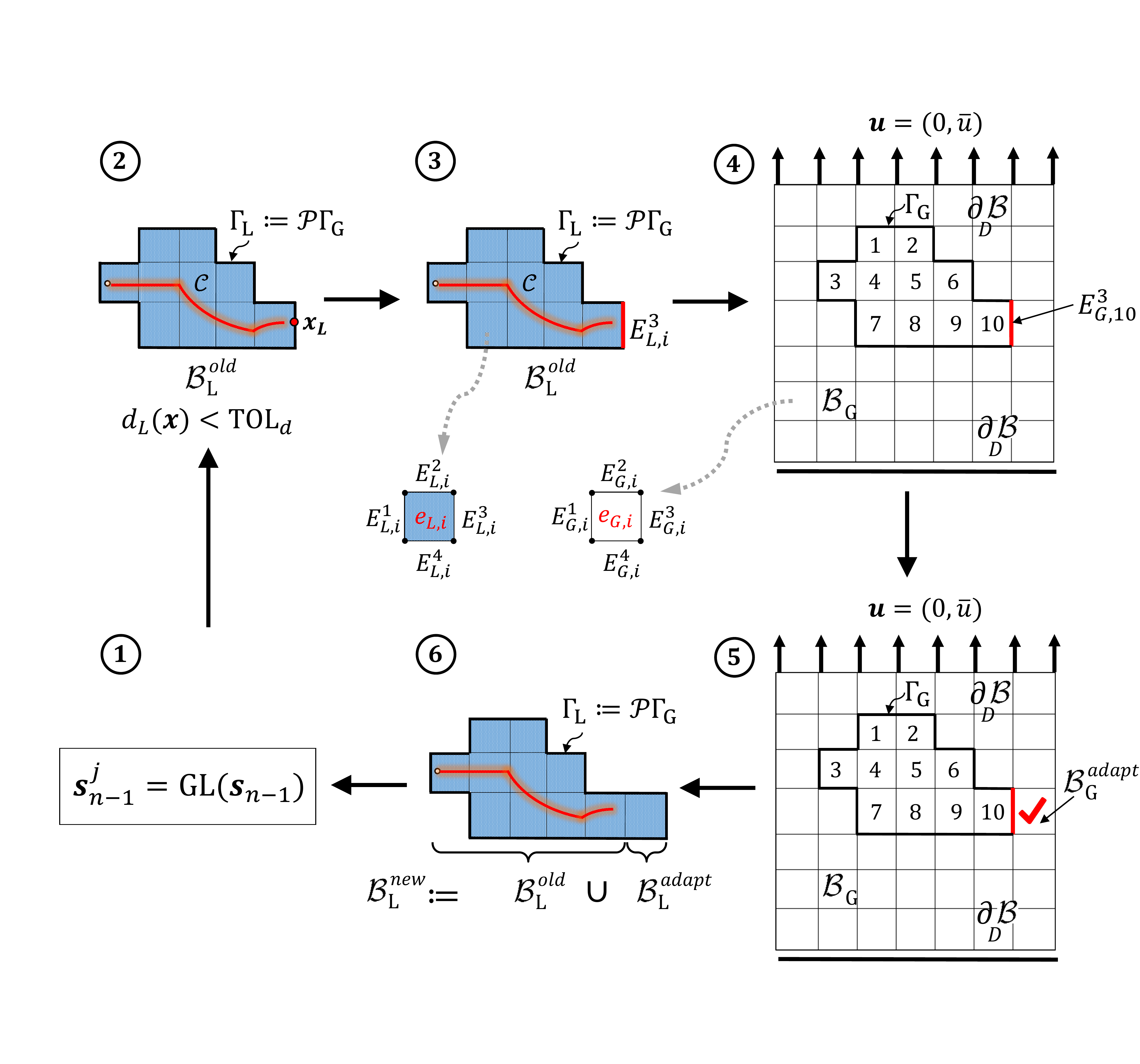}}  
	\caption{{Explanation of the predictor-corrector adaptive scheme, introduced in \cite{NoiiAldakheelWickWriggers2019}.}
	}
	\label{adapt}
\end{figure}

\sectpd[Section3321]{Fluid pressure}
Following the same procedure as above, the Robin-type boundary condition at the global level for the pressure field yields
\begin{equation*}
\int_\Gamma \lambda^{p,k}_C \;\delta {p}_\Gamma \,\mathrm{d}A+\KIA_G^p\int_\Gamma p_\Gamma^{k}\; \delta {\lambda}_L^p \,\mathrm{d}A={\Lambda}^{p,k}_G,
\label{RBC14}
\tag{$\widetilde {\text{C}}_7$}
\end{equation*}
\begin{equation*}
\int_\Gamma (p_\Gamma^{k,\frac{1}{2}}-p^{k}_G) \delta {\lambda}_C^p \,\mathrm{d}A=0,
\label{RBC15}
\tag{$\widetilde {\text{C}}_8$}
\end{equation*}  
with
\begin{equation}
{\Lambda}^{p,k}_G:=\Lambda_L(\lambda^{p,k}_L,p_L^{k};\KIA_G^p)=
\KIA_G^p\int_\Gamma p_L^{k} \delta {\lambda}_L^p \,\mathrm{d}A-\int_\Gamma \lambda^{p,k}_L \delta {p}_\Gamma \,\mathrm{d}A.
\label{RBC16}
\end{equation}
Together with (\ref{GD}), the global system of equations has to be solved for $(\Bvarphi^{k}_G, p_G^k, \Blambda^{\varphi,k}_C, \lambda^{p,k}_C, \Bvarphi_\Gamma^{k}, p_\Gamma^{k})$ for a given $({\BLambda}^{\varphi,k}_G, {\Lambda}^{p,k}_G, \KIA_G^\varphi, \KIA_G^p, \Bvarphi_\Gamma^{k,\frac{1}{2}}, p_\Gamma^{k,\frac{1}{2}})$. Here, ($\KIA_G^\varphi$, $\KIA_G^p$, ${\BLambda}^{\varphi,k}_G$, ${\Lambda}^{p,k}_G$) stand for global Robin-type parameters. 

Based on the new boundary conditions provided in ($\widetilde {\text{C}}_1$) -- ($\widetilde {\text{C}}_8$) the imbalanced quantities in the Global-Local iterations read
\begin{equation}
\int_\Gamma (\Bvarphi_\Gamma^k-\Bvarphi_{\Gamma}^{k,\frac{1}{2}}) \cdot \delta {\Blambda}_L^\varphi \,\mathrm{d}{A} \neq 0 
\AND
\int_\Gamma (p_\Gamma^k - p_\Gamma^{k,\frac{1}{2}}) \delta {\lambda}_L^p \,\mathrm{d}{A} \neq 0 
\label{RBC17}
\end{equation}
For the specific Robin-type boundary conditions, we can resolve (\ref{RBC17}) such that this term does not produce any error in the iterative procedure. To do so, following our previous work \cite{NoiiAldakheelWickWriggers2019}, the global and local augmented stiffness matrices for both, the deformation and pressure fields, within the Robin-type boundary conditions are given by
\begin{equation}\label{eq_aug_stiff}
\fterm{
	\KIA_G={\bm L^T_L}{\bm T^{-T}_L}{\bm {\mathcal{S}}}_L \quad \mbox{and} \quad \KIA_L:={\bm {\mathcal{S}}}_C.
}
\end{equation}
$\KIA_G$ and $\KIA_L$ can be seen as augmented stiffness matrices regularize the Jacobian stiffness matrix at the global and local levels, respectively. For details on the derivation of those matrices, we refer the interested reader to \cite{NoiiAldakheelWickWriggers2019}.

\begin{Remark}
	\label{C-choise}
	Note that the choice of the coupling equations (\ref{Coupl1}) -- (\ref{Coupl6})
	at the local and global level for the Robin-type boundary conditions are the outcome of precise investigation of different combinations. However other choices are also possible, but one needs to adapt/derive the imbalance equations in \req{RBC17} accordingly.
\end{Remark}

The detailed Global-Local formulation using Robin-type boundary 
conditions is depicted in {Algorithm} \ref{TableGL}. The Global-Local setting provides a generic two-scale finite element algorithms that enables capturing local non-linearities.

%
\sectpb[Section34]{Predictor-Corrector mesh adaptivity}
The Global-Local approach is augmented by a {\it dynamic allocation} of a local state using an adaptive scheme which has to be performed at time step $t_{n}$ in {Algorithm} \ref{TableGL}. By the adaptivity procedure, we  mean: ({a}) to
determine which global elements need to be refined; ({b}) to create the new fictitious
and local domain, see Fig. \ref{adapt}; ({c}) to determine a new local interface; ({d}) to interpolate the old global solution. For details regarding the predictor-corrector adaptive scheme applied to the global-local formulation, we refer the interested reader to \cite{Wick15Adapt,NoiiAldakheelWickWriggers2019} and Algorithm 2-3 therein.

\newpage

\begin{algorithm}[H]\small
	\caption{\em Global-Local iterative scheme combined with Robin-type boundary conditions.}
	\label{TableGL}
	{
		\begin{tabular}{l}
			{\bf Input:} loading data $(\bar{\Bvarphi}_{n},\bar{p}_{n})$ on $\partial_D \calB$; \\
			\hspace{1.37cm}solution $\BfrakP_{n-1}:=(\Bvarphi_{G,n-1}, p_{G,n-1}, \Bvarphi_{L,n-1}, p_{L,n-1}, d_{L,n-1}, \Bvarphi_{\Gamma,n-1}, p_{\Gamma,n-1}, \Blambda^\varphi_{C,n-1},\Blambda^\varphi_{L,n-1},$ \\ 
			$\qquad \qquad\qquad\qquad \lambda^p_{C,n-1}, \lambda^p_{L,n-1})$ and $\calH_{L,n-1}$ from step $n-1$. \\ [0.2cm]
			
			{{\bf Global-Local iteration} $k\geq1$}:\\[0.2cm]
			
			\quad \underline{Local boundary value problem}:\\
			
			\quad\quad\;\; \textbullet\; given $\KIA_L^\varphi, \KIA_L^p,{\BLambda}^{\varphi,k-1}_L, {\Lambda}^{p,k-1}_L, \calH_{L,n-1}$; solve \\
			\quad\qquad\;\;${\text{phase-field part:}}$ \;$\displaystyle \int_{\calB_L} \Big[ \Big(2\psi_c \;d_L + 2(d_L-1) \ \calH(\nabla \Bvarphi_L) \Big) \delta d_L  + 2\psi_c \;l^2 \; \nabla d_L \cdot \nabla \delta d_L \Big] dV = 0$,\\
			
			\quad\qquad\;\;${\text{mechanical part:}}$\\
			\quad\qquad
			$\begin{cases} 
			\displaystyle \int_{\calB_{L}} \BP(\nabla \Bvarphi_L, p_L, d_L):\nabla \delta\Bvarphi_L \mathrm{d}V - \int_{\Gamma_{L}} \Blambda_L^\varphi \cdot \delta {{\Bvarphi}}_G\,\mathrm{d}A=0, \\
			\displaystyle \int_\Gamma \Blambda_L^\varphi\cdot \delta {\Bvarphi}_\Gamma \,\mathrm{d}A+\KIA_L^\varphi\int_\Gamma \Bvarphi_\Gamma \cdot \delta {\Blambda}_C^\varphi \,\mathrm{d}A={\BLambda}^{\varphi,k-1}_L \AND  \int_\Gamma (\Bvarphi_\Gamma-\Bvarphi_L) \cdot \delta {\Blambda}_L^\varphi \,\mathrm{d}A=0,		
			\end{cases}$ \\
			
            \quad\qquad\;\;${\text{fluid pressure:}}$\\
\quad\qquad
            $\begin{cases} 
            \displaystyle \int_{\calB_L} \Big[\frac{1}{M}(p_L-p_{L,n}) + B \Big(J(\nabla \Bvarphi_L)-J_n(\nabla \Bvarphi_L)\Big)  \Big]\delta p_L\; \mathrm{d}V\\
            + \int_{\calB_L} \Big[ \Big(\Delta t \;\BK(\nabla \Bvarphi_L, d_L) \;\nabla p_L \Big) \cdot \nabla \delta p_L \Big] \mathrm{d}V
            -\int_{\Gamma_{L}} \lambda_L^p \ \delta {p}_L\,\mathrm{d}A = 0,\\
            \displaystyle \int_\Gamma \lambda_L^p \delta {p}_\Gamma \,\mathrm{d}A+\KIA_L^p\int_\Gamma p_\Gamma \delta {\lambda}_C^p \,\mathrm{d}A={\Lambda}^{p,k-1}_L \AND  \int_\Gamma (p_\Gamma-p_L) \delta {\lambda}_L^p \,\mathrm{d}A=0,		
            \end{cases}$ \\ 			
			
			\quad\quad\quad \;\;{\color{white}\textbullet}\; \;set $(\Bvarphi_L,p_L,d_L,\Bvarphi_\Gamma,p_\Gamma,\Blambda_L^\varphi,\lambda_L^p):=(\Bvarphi_L^k,p_L^k,d_L^k,\Bvarphi_\Gamma^{k,\frac{1}{2}},p_\Gamma^{k,\frac{1}{2}},\Blambda_L^{\varphi,k},\lambda_L^{p,k})$, \\
			
			\quad\qquad \textbullet\; given $({\Bvarphi_L^{k}},p_L^{k},{\Blambda^{\varphi,k}_L},{\lambda^{p,k}_L};\KIA_G^\varphi,\KIA_G^p)$, set\\
			\quad\quad\qquad $\displaystyle {\BLambda}^{\varphi,k}_G
			=\KIA_G^\varphi\int_\Gamma \Bvarphi_L^{k} \cdot \delta {\Blambda}_C^\varphi \,\mathrm{d}A-\int_\Gamma \Blambda^{\varphi,k}_L\cdot \delta {\Bvarphi}_\Gamma \,\mathrm{d}A$\ ; \ $\displaystyle {\Lambda}^{p,k}_G
			=\KIA_G^p\int_\Gamma p_L^{k} \delta {\lambda}_C^p \,\mathrm{d}A-\int_\Gamma \lambda^{p,k}_L \delta {p}_\Gamma \,\mathrm{d}A$. \\[0.2cm]
			
			\quad \underline{Global boundary value problem}:\\
			
			\quad\quad\;\; \textbullet\; given $\KIA_G^\varphi,\KIA_G^p,{\BLambda}^{\varphi,k}_G,{\Lambda}^{p,k}_G,\Bvarphi_\Gamma^{k,\frac{1}{2}},p_\Gamma^{k,\frac{1}{2}}$, solve \\	
			\quad\qquad\;\;${\text{mechanical part:}}$\\
			\hspace{1.cm}
			$\begin{cases} 
			\displaystyle \int_{\calB_G}\BP(\nabla\Bvarphi_G,p_G,0):\nabla\delta {\Bvarphi}_G \,\mathrm{d}V - \int_{\calB_f}\BP(\nabla\Bvarphi_G,p_G,0):\nabla\delta {\Bvarphi}_G \,\mathrm{d}V\\
			-\int_\Gamma \Blambda_C^\varphi \cdot \delta {\Bvarphi}_G \,\mathrm{d}{A}
			-\int_{\Gamma_{N}} \bar{\Btau} \cdot \delta {\Bvarphi}_G \,\mathrm{d}{A}=0,  \\
			\displaystyle \int_\Gamma \Blambda_C^\varphi\cdot \delta {\Bvarphi}_\Gamma \,\mathrm{d}A+\KIA_G^\varphi\int_\Gamma \Bvarphi_\Gamma \cdot \delta {\Blambda}_L^\varphi \,\mathrm{d}A={\BLambda}^{\varphi,k}_G \AND  \int_\Gamma (\Bvarphi_\Gamma^{k,\frac{1}{2}}-\Bvarphi_G) \cdot \delta {\Blambda}_C^\varphi \,\mathrm{d}A=0,	
			\end{cases}$ \\

			\quad\qquad\;\;${\text{fluid pressure:}}$\\
			\hspace{1.cm}
			$\begin{cases} 
			\displaystyle \int_{\calB_G} \Big[\frac{1}{M}(p_G-p_{G,n}) + B \Big(J(\nabla \Bvarphi_G)-J_n(\nabla \Bvarphi_G)\Big)  \Big]\delta p_G\; \mathrm{d}V -\int_{\Gamma} \lambda_C^p \ \delta {p}_G\,\mathrm{d}A\\
			+ \int_{\calB_G} \Big[ (\Delta t \;\BK(\nabla \Bvarphi_G) \nabla p_G) \cdot \nabla \delta p_G \Big] \mathrm{d}V + \int_{\calB_f} \Big[\frac{1}{M}(p_G-p_{G,n}) + B \Big(J(\nabla \Bvarphi_G)\\
			- J_n(\nabla \Bvarphi_G)\Big) 
			\Big]\delta p_G\; \mathrm{d}V + \int_{\calB_f} \Big[ (\Delta t \;\BK(\nabla \Bvarphi_G) \;\nabla p_G) \cdot \nabla \delta p_G -\Delta t \;\bar{r}_F \Big] \mathrm{d}V + \int_{\Gamma_{N}} \bar{f} \;\delta p_G\; dA = 0,\\
			\displaystyle \int_\Gamma \Blambda_C^p \delta {p}_\Gamma \,\mathrm{d}A+\KIA_G^p\int_\Gamma p_\Gamma \delta {\Blambda}_L^p \,\mathrm{d}A={\Lambda}^{p,k}_G \AND \int_\Gamma (p_\Gamma^{k,\frac{1}{2}}-p_G) \delta {\lambda}_C^p \,\mathrm{d}A=0,	
			\end{cases}$ \\
			
			\quad\quad\;\; {\color{white}\textbullet}\; set $(\Bvarphi_G, p_G,\Bvarphi_\Gamma,p_\Gamma,\Blambda_C^\varphi,\lambda_C^p)=:(\Bvarphi_G^k, p_G^k,\Bvarphi_\Gamma^k,p_\Gamma^k,\Blambda_C^{\varphi,k},\lambda_C^{p,k})$, \\
			
			\quad\qquad \textbullet\; given $({\Bvarphi_G^{k}},p_G^k,{\Blambda^{\varphi,k}_C},\lambda^{p,k}_C;\KIA_L^\varphi,\KIA_L^p)$, set \\
			\quad\quad\qquad$\displaystyle {\BLambda}^{\varphi,k}_L
			=\KIA_L^\varphi\int_\Gamma \Bvarphi_G^{k} \cdot \delta {\Blambda}_C^\varphi \,\mathrm{d}A-\int_\Gamma \Blambda^{\varphi,k}_C\cdot \delta {\bm u}_\Gamma \,\mathrm{d}A$ and $\displaystyle {\BLambda}^{p,k}_L
			=\KIA_L^p\int_\Gamma p_G^{k} \delta {\lambda}_C^p \,\mathrm{d}A-\int_\Gamma \lambda^{p,k}_C \delta {p}_\Gamma \,\mathrm{d}A$.\\
			
			\quad\quad\;\; \textbullet\; if fulfilled, set $\BfrakP^k=:\BfrakP_{n}$ and stop; \\
			\quad\quad\quad\;\; {\color{white}\textbullet}\; else $k+1\rightarrow k$. \\[0.2cm]
			
			{\bf Output:} solution $\BfrakP_{n}$ and $\calH_{L,n}$.
			%
		\end{tabular}
	}
\end{algorithm}
%
%
%
\sectpa[Section5]{Numerical Examples}
This section demonstrates the performance of the proposed adaptive
Global-Local approach applied to the phase-field modeling of hydraulic fracture in fluid-saturated porous media. Two numerical model problems for the GL formulations are investigated. A considerable reduction of the computational cost is observed in comparison with the single-scale
solution. The material parameters used in both examples are listed in Table~\ref{material-parameters} and based on \cite{MieheMauthe2015,xia2017phase}. For the numerical simulation all variables in both, the global and local domains, are discretized by bilinear quadrilateral $Q1$ finite elements. The total number of elements for the single-scale problem is 28900 elements and for Global domain is 100 elements. The number of elements for the local domain is determined based on predictor-corrector mesh adaptivity.

\renewcommand{\tablename}{Table}
\setcounter{table}{0}
\begin{table}[t]
	\caption{Material parameters used in the numerical examples based on \cite{MieheMauthe2015,xia2017phase}.}
	\centering
	\begin{tabular}{cclll}
		No.  &Parameter & Name                   & Value    & Unit            \\\hline 
		1.   &$E$        & Young's modulus       & $15.96$    & $\mathrm{GPa}$ \\
		2.   &$\nu$      & Poisson's ratio       & $0.2$    & --                 \\
		3.   &$M$        & Biot's modulus        & $12.5$  & $\mathrm{GPa}$ \\
		4.   &$B$        & Biot's coefficient     & $0.79$   & -- \\
		5.   &$K$& Intrinsic permeability  & $2 \times 10^{-14} $  & $\mathrm{m^2}$ \\
		6.   &$K_c$& Spatial permeability in fracture  & $83.3 $  & $\mathrm{m^3s/kg}$ \\
		7.   &$\zeta$& Permeability transition exponent  & $50 $  & -- \\
		8.   &$\eta_F$  & Dynamic fluid viscosity   & $1 \times 10^{-3}$  & $\mathrm{kg/(m.s)}$ \\
		9.   &$\sigma_c$     & Critical effective stress  & $0.005$ & $\mathrm{GPa}$ \\ \hline
		\label{material-parameters}
	\end{tabular}
\end{table}

\begin{figure}[!b]
	\centering
	{\includegraphics[clip,trim=0cm 14cm 10cm 14cm, width=14.5cm]{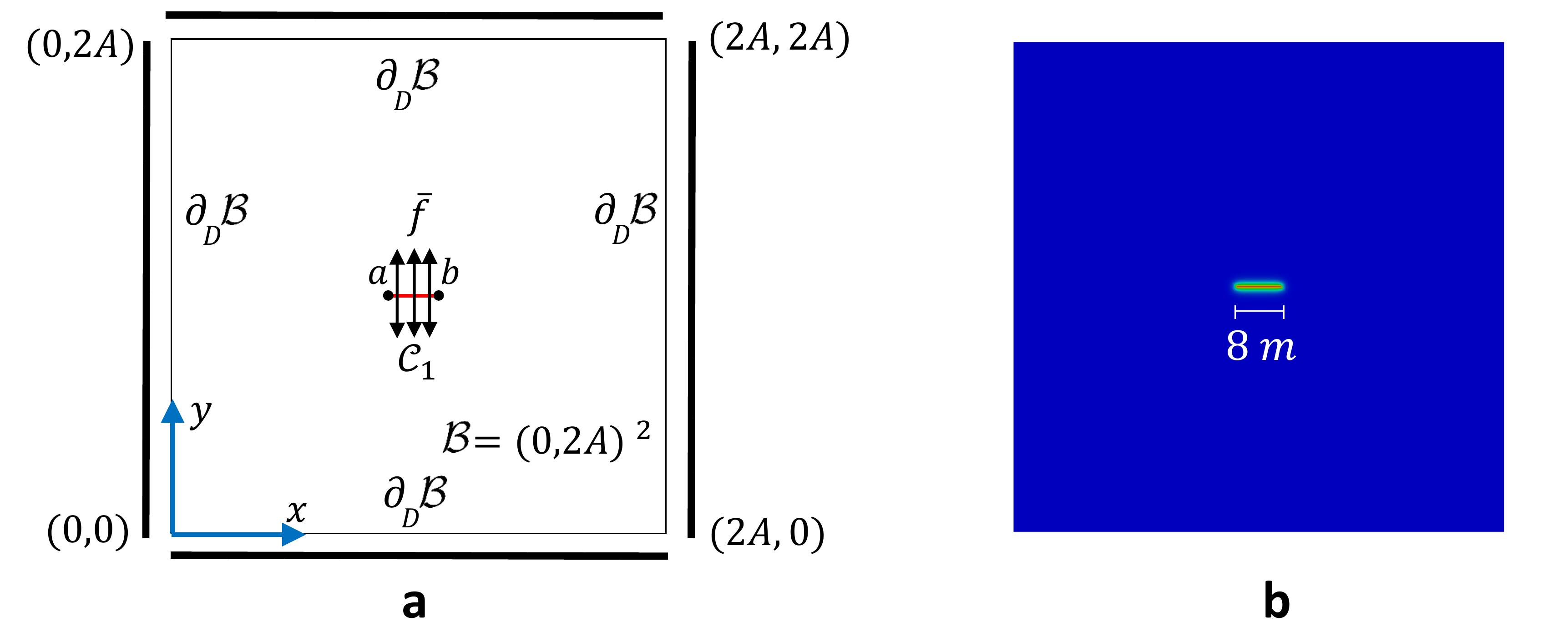}}  
	\caption{Hydraulically induced crack driven by fluid volume injection. (a) Geometry and boundary conditions and (b) described crack phase-field $d$ as a Dirichlet boundary conditions at $t=0~s$.
	}
	\label{example1-a}
\end{figure}

\begin{figure}[!t]
	\centering
	{\includegraphics[width=\textwidth]{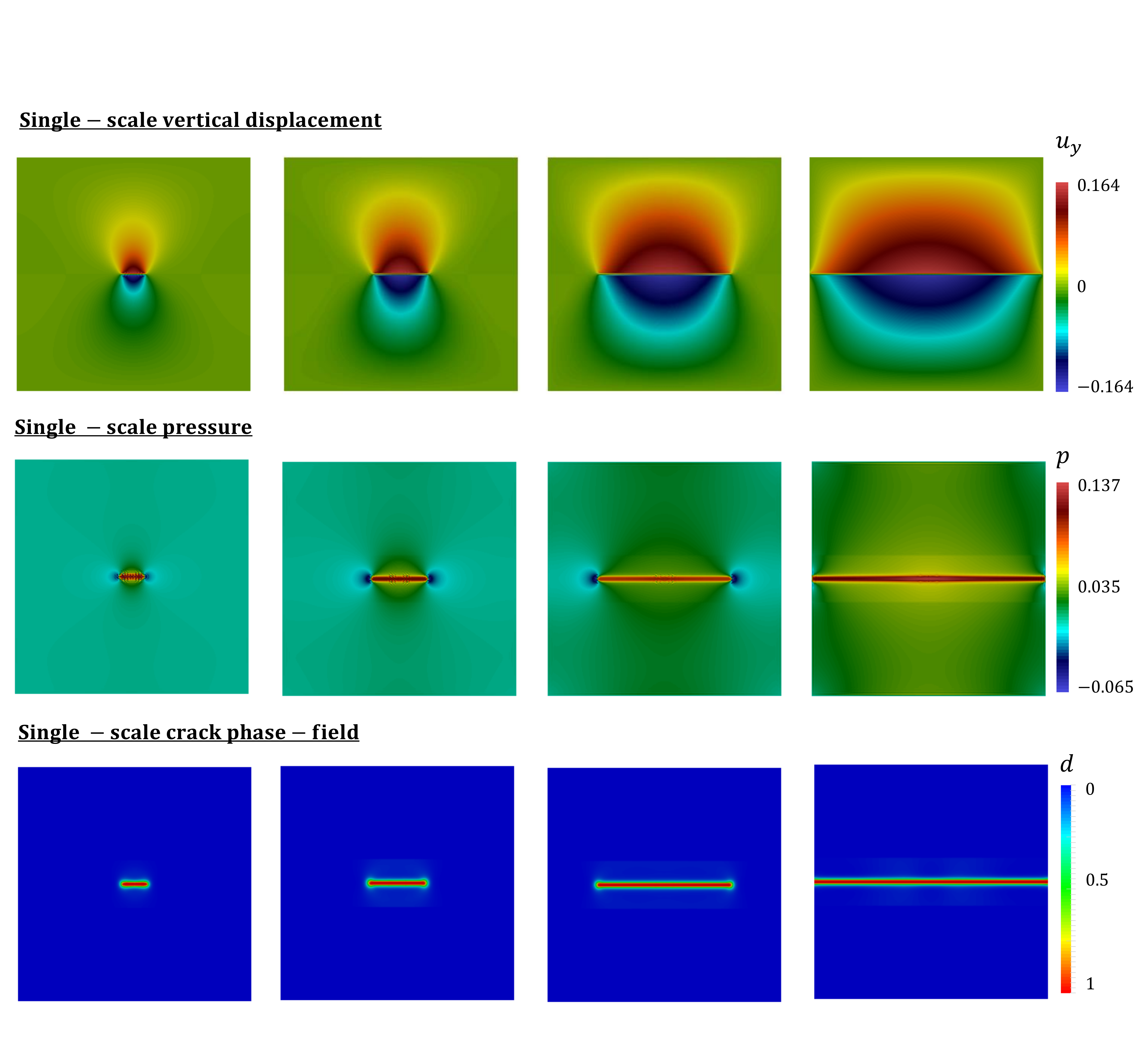}}  
	\caption{Example 1. Single-scale results of the hydraulically induced crack driven by fluid volume injection. Evolution of the vertical displacement $u_y$ (first row), fluid pressure $p$ (second row) and crack phase-field $d$ (third row) for different deformation stages up to final failure at [$t=1.8;6.5;20;48.5~s$].
	}
	\label{example1-b}
\end{figure}

\sectpb{Hydraulically induced crack driven by fluid volume injection}
In the first numerical example, a boundary value problem applied to the square plate is shown in Fig. \ref{example1-a}. We set $A=40~m$ hence $\calB=(0,80)^2~m^2$ that includes a predefined single notch $\calC_1$ of length $8~m$ in the body center with $a =(36,40)$ and $b =(44,40)$, as depicted in Fig. \ref{example1-a}. 
A constant fluid flow of $\bar{f} = 0.002\; m^2/s$ is injected in $\calC_1$. At the boundary $\partial_D\calB$, all the displacements are fixed in both directions and the fluid pressure is set to zero. Fluid injection $\bar{f}$ continues until failure for $T = 49~s$ with time step $\Delta t = 0.1~s$ during the simulation.

We start our analysis by illustrating the single-scale results for different deformations states up to final failure. The vertical displacement $u_y$ (first row), fluid pressure $p$ (second row) and crack phase-field $d$ (third row) evolutions are demonstrated in Figure \ref{example1-b} for four time steps [$t=1.8;6.5;20;48.5~s$]. The crack initiates at the notch-tips due to fluid pressure increase until a threshold energy $\psi_c$ is reached. Thereafter, the crack propagates horizontally in two direction towards the boundaries. In the fractured zone, $p$ is almost constant due to the increased permeability inside the crack. Whereas, low fluid pressure in the surrounding is observed due the chosen small time-step in comparison with the permeability of the porous medium, as outlined in \cite{MieheMauthe2015}. The fluid pressure drops down while the crack propagates further as shown in Figure \ref{example1-b} (second row, middle states). Then, $p$ increases again due to the prescribed fixed boundary conditions $\partial_D\calB$, see Fig. \ref{example1-b} (second row, last state).

\begin{figure}[!t]
	\centering
	{\includegraphics[width=\textwidth]{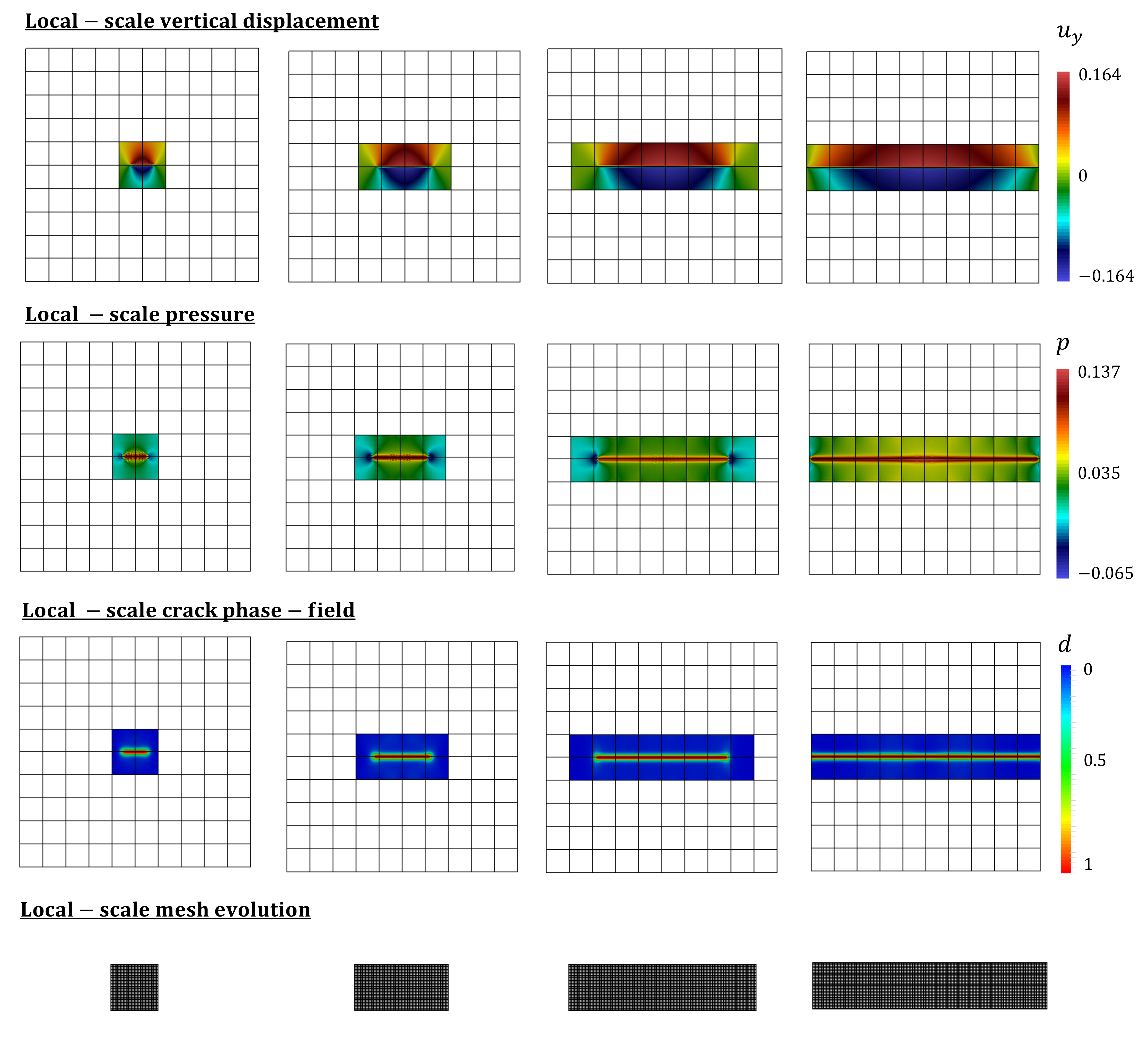}}  
	\caption{Example 1. Local-scale results of the hydraulically induced crack driven by fluid volume injection. Evolution of the vertical displacement $u_y$ (first row), fluid pressure $p$ (second row), crack phase-field $d$ (third row) and local domain (fourth row) for different fluid injection time steps at [$t=1.8;6.5;20;48.5~s$] up to final failure.}
	\label{example1-b-L}
\end{figure}

\begin{figure}[!t]
	\centering
	{\includegraphics[width=\textwidth]{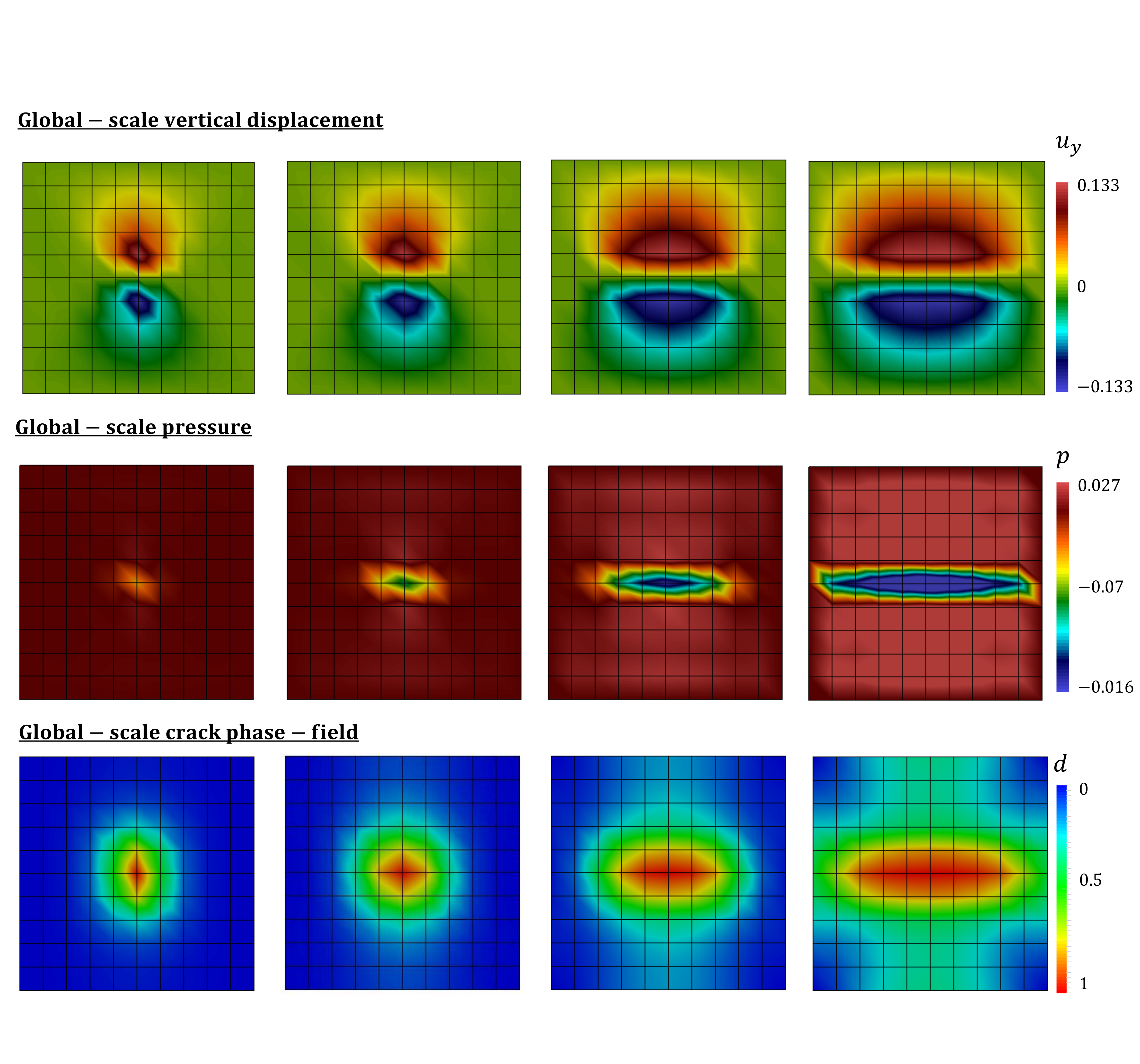}}  
	\caption{Example 1. Global-scale results of the hydraulically induced crack driven by fluid volume injection. Evolution of the vertical displacement $u_y$ (first row), fluid pressure $p$ (second row) and crack phase-field $d$ (third row) for different deformation stages up to final failure at $t=1.8;6.5;20;48.5~s$.
	}
	\label{example1-b-G}
\end{figure}

Next the performance of the Global-Local approach is investigated. To this end, the evolution of the vertical displacement $u_y$ (first row), fluid pressure $p$ (second row), crack phase-field $d$ (third row) and local domain (fourth row) for different fluid injection time steps $t=1.8;6.5;20;48.5~s$ are demonstrated in Fig. \ref{example1-b-L} for the local-scale and in Fig. \ref{example1-b-G} for the homogenized global scale. Hereby, even with less number of elements at the global domain the overall response is qualitatively in a good agreement with the single-scale domain.

\begin{figure}[!t]
\centering
{\includegraphics[clip,trim=0cm 10cm 0cm 16cm, width=16cm]{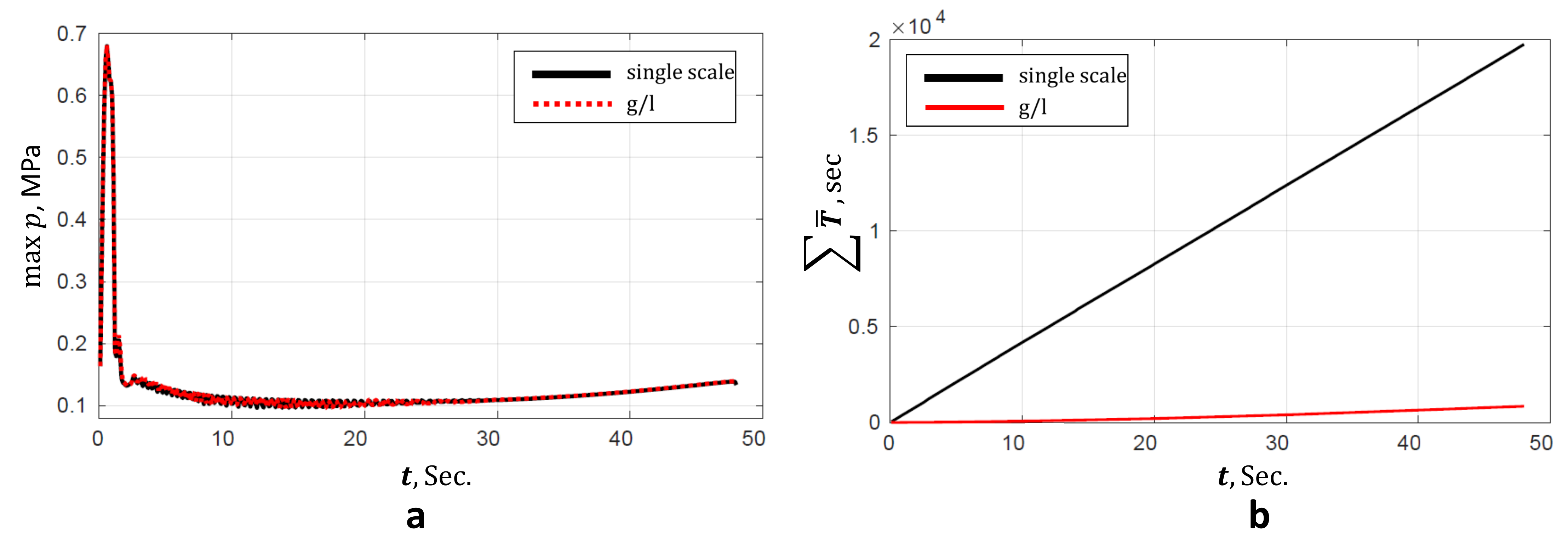}}  
\caption{Example 1. Hydraulically induced crack driven by fluid volume injection. (a) Fluid pressure $p$ within the crack region versus injection time; and (b) computational time-injection time curves in terms of the accumulated time.}
\label{example1-fig1}
\end{figure}

Figure \ref{example1-fig1}a describes the maximum injected fluid pressure within the crack region versus the fluid injection time. For a comparison purpose, the Global-Local and single-scale solutions are both provided. The results obtained from Global-Local formulation are in a good agreement with the single-scale solution. Furthermore, it is noted that the injected fluid pressure increases within the crack region before it reaches to the peak point. Thereafter, as expected a drop of the fluid pressure is observed.

\begin{figure}[!t]
	\centering
	{\includegraphics[clip,trim=0cm 10cm 0cm 16cm, width=16cm]{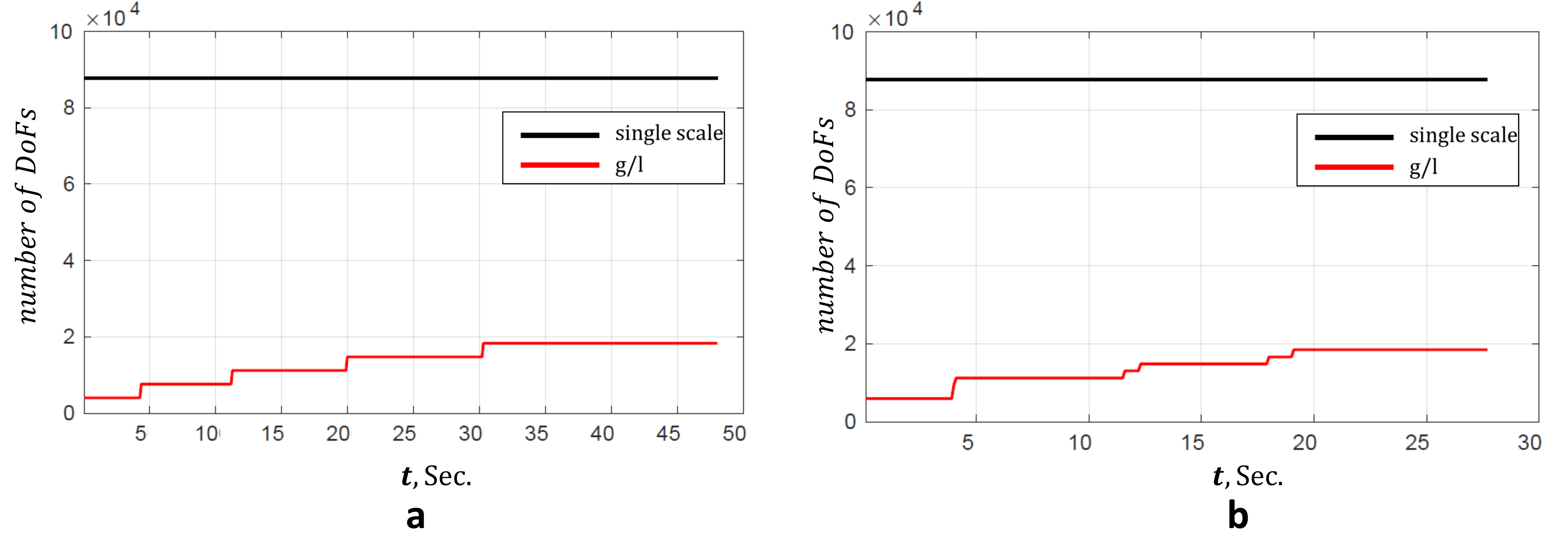}}  
	\caption{Number of degrees of freedom for the single scale problem and Global-Local formulation. (a) Example 1; and (b) Example 2.
	}
	\label{example1-fig2}
\end{figure}

To illustrate the efficiency of the predictor-corrector adaptive scheme, we plot in \ref{example1-fig1}b the corresponding accumulative computational time and in Fig. \ref{example1-fig2}a the total number of unknowns (local and global problems) versus the fluid injection time and compared with the single-scale problem. It can be observed that the total accumulated time for the Global-Local formulations took $849~s$ whereas the single-scale simulation took $19752~s$. Hence, Global-Local formulations performs $23.3$ times faster. 

Furthermore, for each jump in Fig. \ref{example1-fig2}a, the predictor-corrector adaptive scheme is active and applied on the Global-Local scheme which increases the number of degrees of freedoms.
At the complete failure state, i.e. $t=48.5~s$, the total number of local nodes, elements and the degrees of freedoms are 5985, 5780 and 17955, respectively for the Global-local formulations. Whereas for the single-scale the number of nodes, elements and the degrees of freedoms are 29241, 28900 and 87723, respectively.
Hence the Global-Local approach requires significantly less degrees of freedom, as shown in Fig. \ref{example1-fig2}a.

\begin{figure}[!t]
	\centering
	{\includegraphics[clip,trim=0cm 13.5cm 10.5cm 14cm, width=14.5cm]{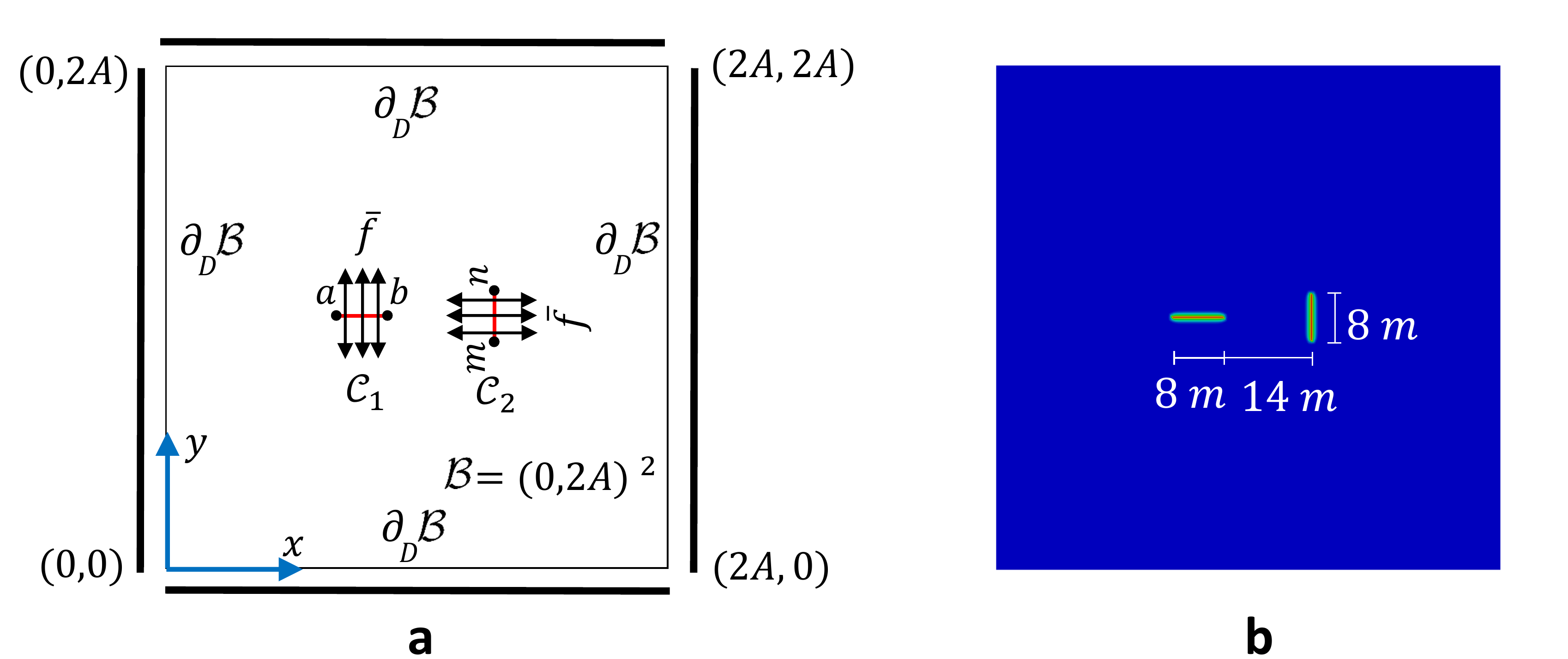}}  
	\caption{Joining of two cracks driven by fluid volume injection. (a) Geometry and boundary conditions; and (b) described crack phase-field $d$ as a Dirichlet boundary conditions at $t=0~s$.
	}
	\label{example2-a}
\end{figure}

\begin{figure}[!t]
	\centering
	{\includegraphics[width=\textwidth]{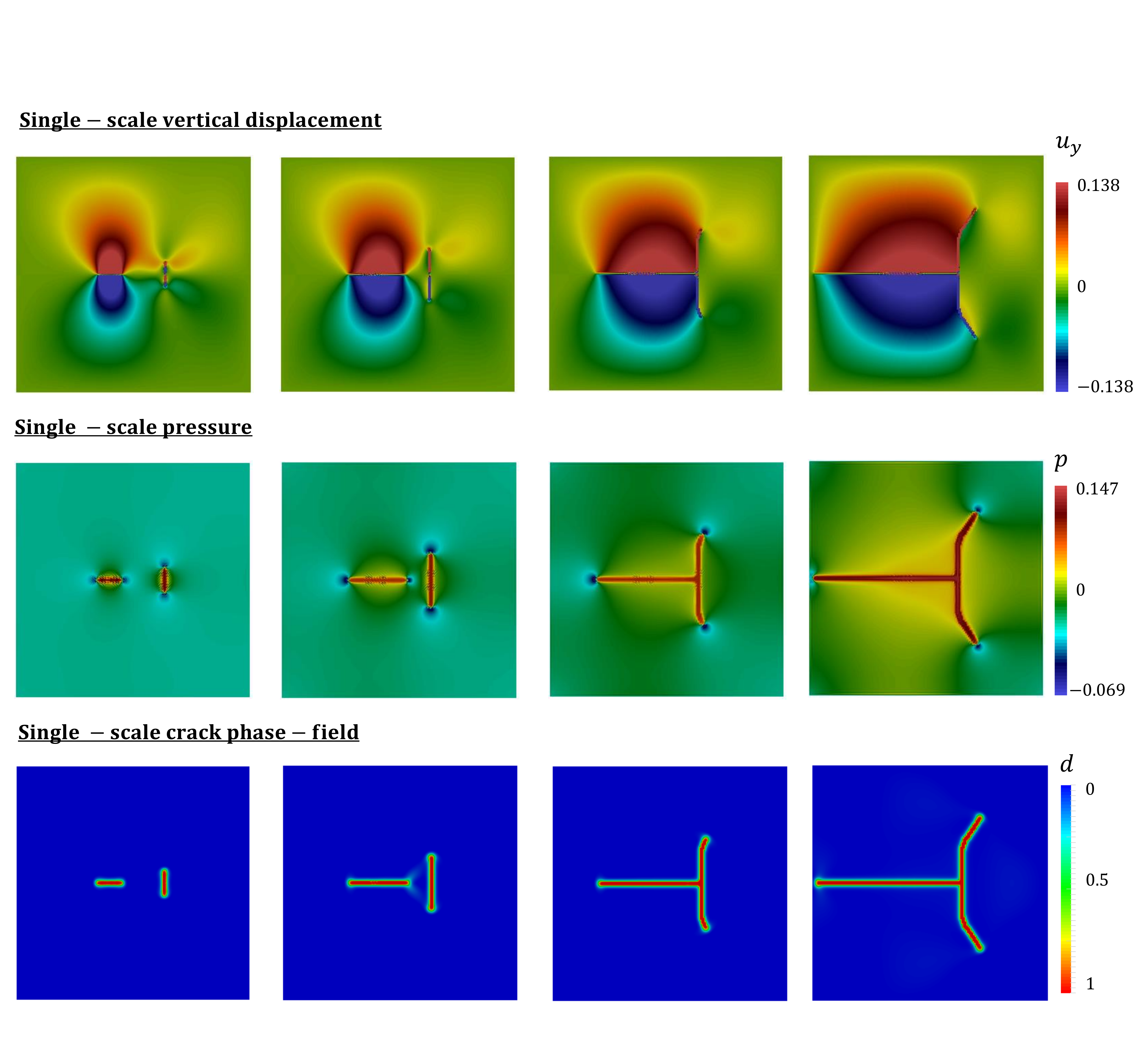}}  
	\caption{Example 2. Single-scale results of the joining cracks driven by fluid volume injection. Evolution of the vertical displacement $u_y$ (first row), fluid pressure $p$ (second row) and crack phase-field $d$ (third row) for different deformation stages up to final failure at [$t=1.8;6.5;13.5;27.7~s$].
	}
	\label{example2-b}
\end{figure}

\begin{figure}[!t]
	\centering
	{\includegraphics[width=\textwidth]{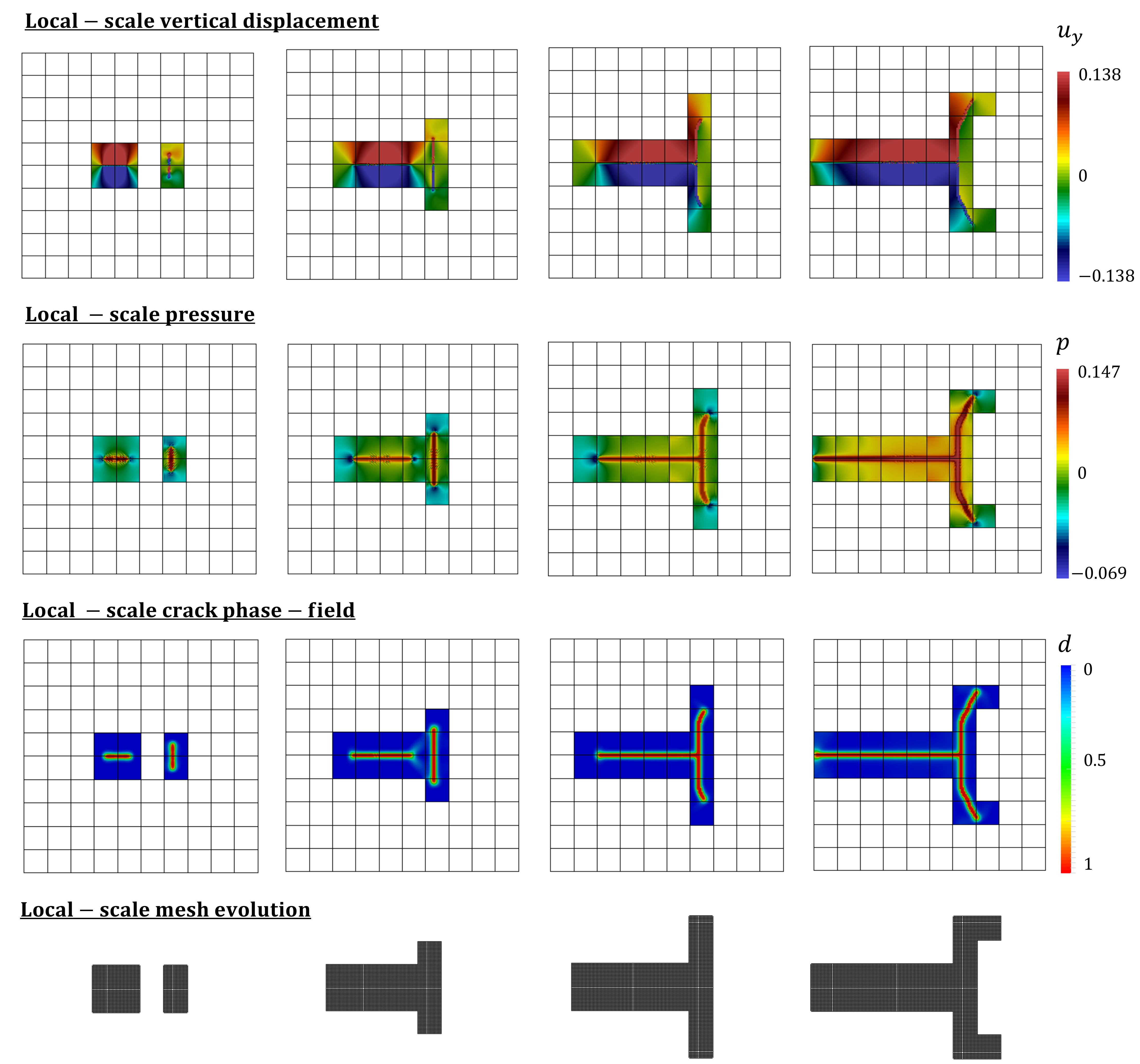}}  
	\caption{Example 2. Local-scale results of the joining cracks driven by fluid volume injection. Evolution of the vertical displacement $u_y$ (first row), fluid pressure $p$ (second row), crack phase-field $d$ (third row) and local domain (fourth row) for different deformation stages up to final failure at [$t=1.8;6.5;13.5;27.7~s$].
	}
	\label{example2-bL}
\end{figure}

\begin{figure}[!t]
	\centering
	{\includegraphics[width=\textwidth]{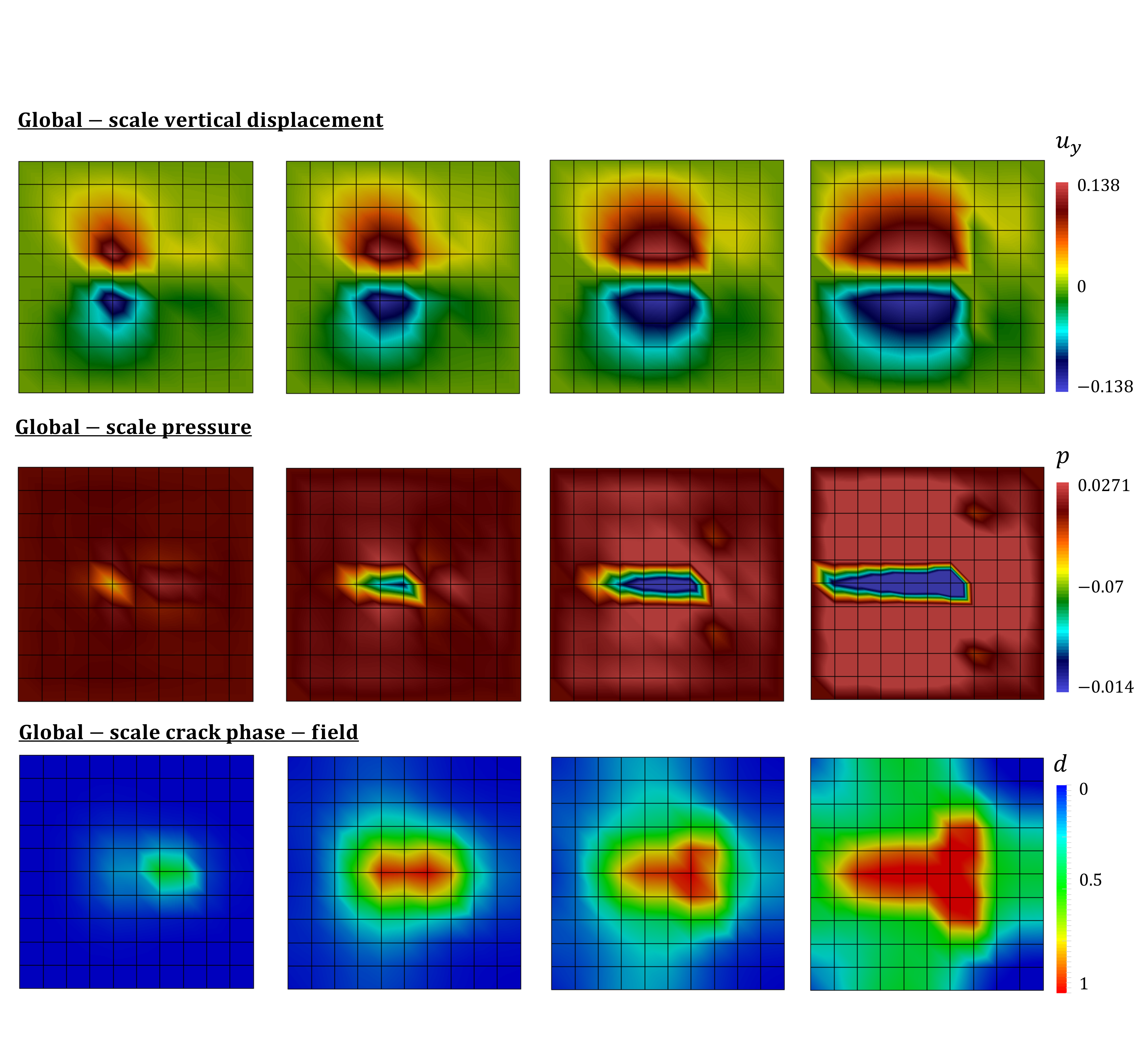}}  
	\caption{Example 2. Global-scale results of the joining cracks driven by fluid volume injection. Evolution of the vertical displacement $u_y$ (first row), fluid pressure $p$ (second row) and crack phase-field $d$ (third row) for different fluid injection stages up to final failure at [$t=1.8;6.5;13.5;27.7~s$].
	}
	\label{example2-bG}
\end{figure}

\begin{figure}[!t]
	\centering
	{\includegraphics[clip,trim=0cm 10cm 0cm 16cm, width=16cm]{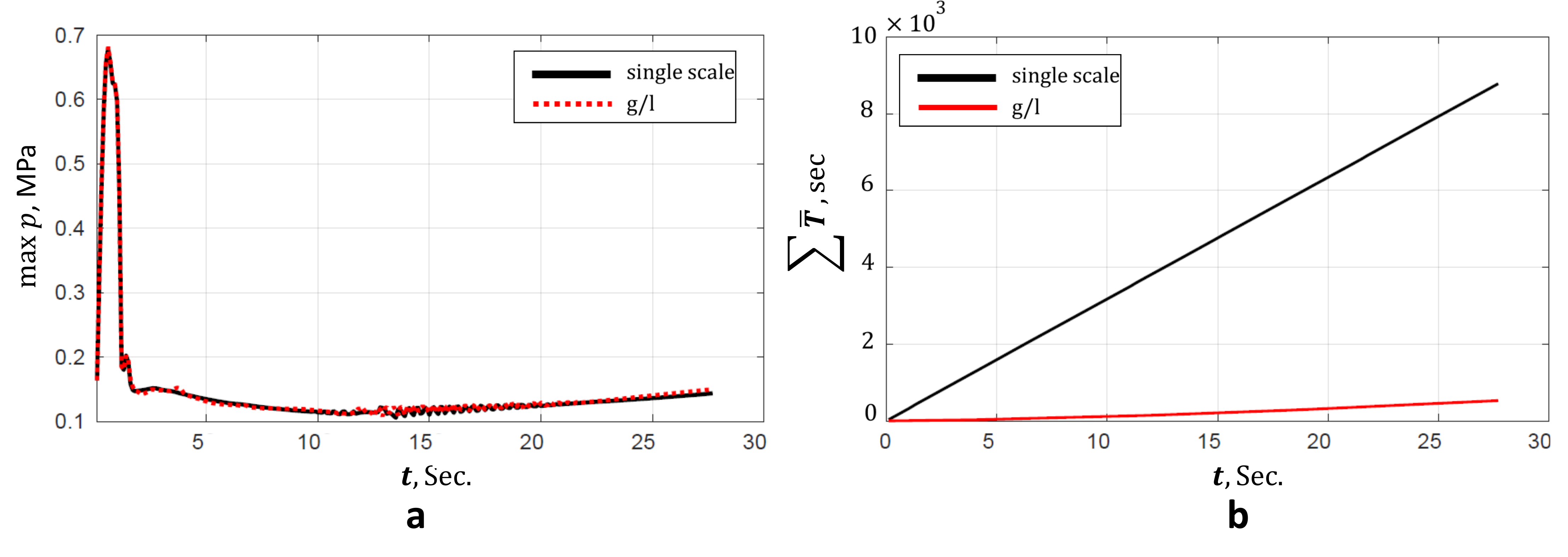}}  
	\caption{Example 2. Joining of two cracks driven by fluid volume injection. (a) Fluid pressure $p$ within the crack region versus injection time; and (b) computational time-injection time curves in terms of the accumulated time.
	}
	\label{example2-fig1}
\end{figure}

\sectpb{Joining of two cracks driven by fluid volume injection}
The second example is concerned with the capability of the proposed GL approach for handling \textit{coalescence} and \textit{merging} of \textit{crack paths} in the local domains. 
Crack-initiation and curved-crack-propagation, representing a mixed-mode fracture, are predicted with a Global-Local formulation. 

The geometrical setup and the loading conditions of the specimen is similar to the benchmark problem of \cite{WickLagrange2014} and depicted in Fig. \ref{example2-a}. We keep all parameters and loading as in the previous example. The first crack $\calC_1$ is located near the middle of the domain with coordinates $a =(28,40)$ and $b =(36,40)$. The second crack $\calC_2$ is vertically-oriented at $n=(50,44)$ and $m=(50,36)$ with a distance of $14\;m$ from $\calC_1$. A constant fluid flow of $\bar{f} = 0.002\; m^2/s$ is injected in $\calC_1$ and $\calC_2$ as sketched in Fig. \ref{example2-a}. At the boundary $\partial_D\calB$, all the displacements are fixed in both directions and the fluid pressure is set to zero. Fluid injection $\bar{f}$ continues until failure for $T = 28~s$ with time step $\Delta t = 0.1~s$ during the simulation.

Figure \ref{example2-b} shows the evolutions of the fluid pressure $p$ (first row) and the crack phase-field $d$ (second row) for the single-scale problem at different times  [$t=1.8;6.5;13.5;27.7~s$]. Here the crack propagates from the notches. We again observe nearly constant fluid pressure in the fractured area ($d=1$), whereas outside the crack zone $p$ is much lower, see \ref{example2-a} (first row).

The local-scale results with the corresponding mesh are depicted in Fig. \ref{example2-bL} for different fluid injection stages. Hereby, the vertical displacement $u_y$ (first row), fluid pressure $p$ (second row), crack phase-field $d$ (third row) and local domain (fourth row) evolutions of the Global-Local formulation are demonstrated in Fig. \ref{example1-b-L} for four time steps [$t=1.8;6.5;13.5;27.7~s$]. It is remarkably observed that the Global-Local approach augmented with predictor-corrector mesh adaptivity leads to the optimum number of elements to be used for the simulation, hence reducing additional cost. Additionally, note that extending the reservoir domain will significantly increase the computational cost for the single-scale problem (due to increase the number of elements) but this will not change the computational cost for Global-Local formulation, thus applicable for the real large structure. Therefore, localize effect (crack phase-field) which increase the computational cost is only considered within local domain and hence globally reduce the computational time. Another advantage of using the GL formulation is its capability of capturing the crack initiation and propagation at the homogenized global scale even with less number of elements as illustrated in Fig. \ref{example2-bG}.

Next, the maximum injected fluid pressure within the crack region is analyzed versus the fluid injection time in Figure \ref{example2-fig1}a. The results obtained from Global-Local formulation are in a good agreement with the single scale solution. Figure \ref{example2-fig1}b represents the corresponding accumulative computational time (i.e. CPU simulation time), per injection fluid time. In this study, we observed that the total accumulated time for the GL approach {took} $529~s$ whereas the single-scale problem took $8784~s$. Hence, Global-Local formulations performs $16.6$ times faster.

Finally, Fig. \ref{example1-fig2}b demonstrates the total number of degrees of freedoms versus the fluid injection time for GL scheme and the single-scale problem. At the complete failure state, i.e. $t=27.7 \;sec$, the total number of local nodes, elements and the degrees of freedoms for the GL method are 6036, 5780 and 18108, respectively. Whereas for the single-scale formulation the number of nodes, elements and the degrees of freedoms are 29241, 28900 and 87723, respectively. Thus the Global-Local approach requires significantly less degrees of freedom.

\sectpa[Section7]{Conclusion}
In this work, we developed a global-local approach for pressurized fractures in porous media. This approach has the potential to tackle practical field problems in which a large reservoir might be considered and fracture propagation is a localized phenomenum.

First, we discussed the governing equations in Section 2. Then, we developed 
the main algorithm in Section 3. Therein, the coupling between the local and the global domain is formulated with Robin-type interface conditions. In our numerical tests, we have shown that the GL approach besides its feasibility for having two different finite element models for the global and local domain, enabled computations with legacy codes. Additionally, it required significantly less degrees of freedom than the single-scale formulation, leads to a remarkable reduction of the computational time. In this regard, the GL approach was $23$ times faster than the standard phase-field formulation (single-scale solution) in the first example and $16$ times faster in the second numerical test. Yet, an excellent performance of the proposed framework in all examples was observed.

\subsection*{Acknowledgment}
F. Aldakheel was founded by the \textsc{Priority Program DFG - SPP 2020} within its second funding phase. N. Noii was partially supported by the \textsc{Priority Program DFG - SPP 1748}. T. Wick and P. Wriggers were funded by the Deutsche Forschungsgemeinschaft (DFG, German Research Foundation)
under Germany's Excellence Strategy within the \textsc{Cluster of
Excellence PhoenixD (EXC 2122)}.

\bibliographystyle{bibls1}
\bibliography{./lit.bib}

\end{document}